\documentclass[11pt,twoside]{article}
\usepackage{amsmath, amssymb, amsthm}
\usepackage{graphicx} 

\theoremstyle{plain}

\newtheorem{thm}{Theorem}[section]
\newtheorem{prop}[thm]{Proposition}
\newtheorem{lemma}[thm]{Lemma}

\newtheorem{corollary}[thm]{Corollary}

\theoremstyle{definition}

\begin{document}


\def\Cal#1{{\cal#1}}
\def\<{\langle}\def\>{\rangle}
\def\what{\widehat}\def\wtil{\widetilde}
\def\Z{{\mathbb Z}}\def\N{{\mathbb N}} \def\C{{\mathbb C}}
\def\Q{{\mathbb Q}}\def\R{{\mathbb R}} \def\E{{\mathbb E}}
\def\P{{\mathbb P}}\def\T{{{\rm Trans}}}\def\Conj{{{\rm Conj}}}
\def\A{{\cal A}}

\let\liff\Longleftrightarrow \let\imply\Rightarrow

\def\Proof{\paragraph{Proof.}}
\def\Remark{\paragraph{Remark.}}
\def\endproof{\hfill$\square$\break\medskip}
\def\noproof{\hfill$\square$\break}
\def\noi{\noindent}

\def\notation{\paragraph{Notation.}}
\def\ackn{\paragraph{Acknowledgement.}}

\let\demph\textbf
\def\span{\text{span}}
\def\CE{{\cal E}}
\def\Aut{\text{\sl Aut}}\def\Out{\text{\sl Out}}\def\Inn{\text{\sl Inn}}
\def\id{\text{id}}
\def\ov{\overline}
\def\length{\text{\sl length}}
\def\im{\text{\sl Im}}

\def\al{\alpha}                 \def\be{\beta}
\def\ga{\gamma}			\def\Ga{\Gamma}
\def\sig{\sigma}\def\Sig{\Sigma}			\def\de{\delta}
\def\ep{\epsilon}               \def\varep{\varepsilon}
%

\def\Bn{{{\cal B}_n}}
\def\Sn{{{\cal S}_n}}

\title{\bf{Artin groups of type $B$ and $D$}}
 
\author{
\textsc{John Crisp and Luis Paris}\\
\\
\emph{Laboratoire de Topologie}\\
\emph{Universit\'e de Bourgogne}\\
\emph{UMR 5584 du CNRS}\\
\emph{B.P. 47870, 21078 Dijon cedex}\\
\emph{FRANCE}\\
\\
\texttt{jcrisp@u-bourgogne.fr}\\
\texttt{lparis@u-bourgogne.fr}
}

\date{\emph{(\today )}}

\maketitle

\begin{abstract} 
We show that each of the Artin groups of type $B_n$ and $D_n$ can be presented as a semidirect 
product $F \rtimes \Bn$, where $F$ is a free group and $\Bn$ is the $n$-string braid group. We 
explain how these semidirect product structures arise quite naturally from fibrations, and observe that, 
in each case, the action of the braid group $\Bn$ on the free group $F$ is classical.
We prove that, for each of the semidirect products, the group of automorphisms which leave 
invariant the normal subgroup $F$ is small: namely, $\Out(A(B_n),F)$ has order $2$, and 
$\Out(A(D_n),F)$ has order 4 if $n$ is even and 2 if $n$ is odd.
It is known that the Artin group of type $D_n$ may be viewed as an index 2 subgroup of the $n$-string 
braid group over some orbifold. Applying the same techniques, we show that this latter group 
has an outer automorphism group of order 2.
Finally, we determine the automorphism groups of all Artin groups or rank 2.
\end{abstract}



\section{Introduction}\label{sect0}

\begin{figure}[ht]
\begin{center}
\includegraphics[width=15cm]{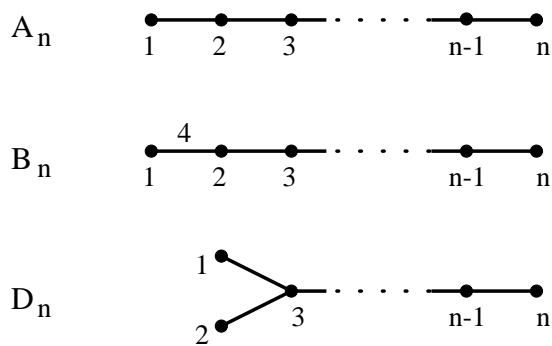}
\end{center}
\caption{Coxeter graphs of type $A_n$, $B_n$ and $D_n$.}\label{fig1}
\end{figure}

Let $S$ be a finite set. A {\it Coxeter matrix} over $S$ is a matrix 
$M=(m_{\alpha\,\beta})_{\alpha, \beta \in S}$ indexed by the elements of $S$
such that $m_{\alpha\,\alpha}=1$ for all $\alpha\in 
S$, and $m_{\alpha\, \beta} = m_{\beta\, \alpha} \in \{2,3,4, \dots, +\infty\}$ for all $\alpha, 
\beta \in S$, $\alpha \neq \beta$. A Coxeter matrix $M=(m_{\alpha\, \beta})$ is usually 
represented by its {\it Coxeter graph}, $\Gamma$. This (labelled) graph is defined by the 
following data. $S$ is the set of vertices of $\Gamma$. Two vertices $\alpha, \beta$ are joined 
by an edge if $m_{\alpha\, \beta} \ge 3$, and this edge is labelled by $m_{\alpha\, \beta}$ if 
$m_{\alpha\,\beta} \ge 4$. For $\alpha, \beta \in S$ and $m\in \Z_{\ge 2}$ we denote by 
$w(\alpha, \beta: m)$ the word $\alpha \beta \alpha \dots$ of length $m$. Define the {\it Artin 
group} of type $\Gamma$ to be the (abstract) group $A(\Gamma)$ presented by
$$
A(\Gamma) = \< S\ |\ w(\alpha, \beta: m_{\alpha\, \beta}) = w(\beta, \alpha: m_{\alpha\, \beta})\ 
{\rm for}\ \alpha \neq \beta\ {\rm and}\ m_{\alpha\, \beta} <+\infty \>\,.
$$
The {\it Coxeter group} $W(\Gamma)$ of type $\Gamma$ is the quotient of $A(\Gamma)$ by the 
relations $\alpha^2=1$, $\alpha \in S$. The cardinal of $S$ is called the {\it rank} of 
$A(\Gamma)$. We say that $A(\Gamma)$ is of {\it spherical type} if $W(\Gamma)$ is finite, and 
that $A(\Gamma)$ is {\it irreducible} if $\Gamma$ is connected. Note that, if $\Gamma_1, 
\Gamma_2, \dots, \Gamma_\ell$ are the connected components of $\Gamma$, then $A(\Gamma)= 
A(\Gamma_1) \times A(\Gamma_2) \times \dots \times A(\Gamma_\ell)$.
If $\Gamma$ is the graph $A_{n-1}$ shown in Figure \ref{fig1}, then $A(\Gamma)=A(A_{n-1})$ 
is the Artin braid group on $n$ strings, which we shall also denote by $\Bn$, and $S=\{ \alpha_1, 
\dots, \alpha_{n-1}\}$ is the set of standard generators of $\Bn$. Moreover, the group $W(A_{n-1})$ is the 
symmetric group on $\{1, \dots, n\}$, which we shall also denote by $\Sn$, and the $\alpha_i$ 
correspond in this group to the transpositions $(i,i+1)$, $i=1, \dots, n-1$.

Artin groups, also called {\it generalized braid groups}, were first introduced by Tits 
\cite{Tit} as extensions of Coxeter groups. Later, Brieskorn \cite{Br1} gave a topological 
interpretation of the Artin groups of spherical type in terms of regular orbit spaces as 
follows. 
Define a {\it finite reflection group} to be a finite subgroup $W$ of $O(n,\R)$ generated by 
reflections, where $n$ is some positive integer. A classical result due to Coxeter \cite{Cox1}, 
\cite{Cox2} states that $W$ is a finite reflection group if and only if $W$ is a finite Coxeter group. 
Assume this is the case. Let $\A (\Gamma)$ be the set of reflecting hyperplanes of $W=W(\Gamma)$, 
and, for $H \in \A (\Gamma)$, let $H_\C$ denote the hyperplane of $\C^n$ having the same equation 
as $H$. Let
$$
M(\Gamma)= \C^n \setminus \left( \bigcup_{H \in \A (\Gamma)} H_\C \right)\,.
$$
Then $M(\Gamma)$ is a connected submanifold of $\C^n$, the group $W(\Gamma)$ acts freely on 
$M(\Gamma)$, and the quotient $N(\Gamma)= M(\Gamma) / W(\Gamma)$ is isomorphic to the complement 
in $\C^n$ of an algebraic variety called {\it discriminantal variety} of type $\Gamma$. Now, by a 
theorem of Brieskorn \cite{Br1}, the fundamental group of $N(\Gamma)$ is isomorphic to the Artin 
group $A(\Gamma)$. In the case $\Gamma=A_{n-1}$, $W(A_{n-1})= \Sn$ acts on $\R^n$ by permuting 
the coordinates, the reflections in $\Sn$ are exactly the transpositions, and $N(A_{n-1})$ is the space 
of configurations of $n$ (unordered) points in $\C$, whose fundamental group is well-known to be 
the braid group $\Bn$.

Since the work of Brieskorn and Saito \cite{BrSa} and that of Deligne \cite{Del}, the 
combinatorial theory of Artin groups of spherical type has been well studied. 
In particular, these groups  are known to be biautomatic (see \cite{Ch1}, \cite{Ch2}).
The finite irreducible Coxeter groups, and therefore the irreducible Artin groups of spherical 
type, were classified by Coxeter \cite{Cox2}. These consist of the three infinite families of Coxeter groups 
defined by the Coxeter graphs $A_n$, $B_n$ and $D_n$ of Figure \ref{fig1}, the dihedral groups
of order $2m$ for $m\geq 3$, associated to the Artin groups of rank 2 (with
$m_{\alpha\, \beta} = m_{\beta\, \alpha} = m \neq +\infty$), as well as the 6 so-called sporadic 
reflection groups. 

In this paper,
for each of the Artin groups of type $B_n$ and $D_n$ we present semidirect
product structures $F\rtimes\Bn$, where $F$ is a free group 
and $\Bn$ is the $n$-string braid group.
We explain how these structures arise quite naturally from fibrations of the 
regular orbits spaces $N(B_n)$ and $N(D_n)$ based on $N(A_{n-1})$. 
Moreover, we observe that, in each case, the action of the braid group $\Bn$ on the free group $F$ is
classical:
for type $B_n$ it is Artin's representation \cite{Ar1,Ar2}, and for type $D_n$ it comes from the monodromy action
on the Milnor fibre of the singularity of type $A_{n-1}$ (see Perron-Vannier \cite{PV}).

The automorphism groups of the (spherical type) Artin groups are just beginning to be explored. Artin's 
1947 paper \cite{Art3} was motivated by the problem of determining the automorphism groups of the 
braid groups (it is explicit in the introduction). However, the problem itself was only solved 37 
years later by Dyer and Grossman \cite{DG}  who proved that the outer automorphism group of the 
braid group $\Bn$ is of order 2 generated by the automorphism which sends each standard generator 
to its inverse. Until now, except for the braid groups, the only known significant result on the 
automorphism groups of (spherical type) Artin groups is an extension of Artin's results of 
\cite{Art3} to all irreducible Artin groups of spherical type (see \cite{CoPa}).
In this paper we prove that, for each of the above semidirect products,
the group of automorphisms which leave invariant the normal subgroup $F$ is small:
namely $\Out(A(B_n),F)$ has order $2$
(see Theorem \ref{Bautos}), and $\Out(A(D_n),F)$ has order $4$ if $n$ is even, and $2$ if $n$ is odd 
(see Theorem \ref{AutKFDn}). It is an open question as to whether, in either the $B_n$ or the $D_n$ case,
the free subgroup $F$ is a characteristic subgroup of the Artin group.

The Artin group of type $D_n$ may be viewed as an index 2 subgroup 
of the $n$-string braid group over an orbifold (namely the complex plane
with a single orbifold point of degree 2, see Allcock \cite{All}). 
This latter group can be presented as a semidirect product $K \rtimes \Bn$, where $K$ is the free 
product of $n$ copies of $C_2= \{\pm 1\}$, and it shall play a major role in our study of the 
Artin groups of type $D_n$. We also show in the paper that the outer automorphism group of this 
group is of order 2.

We determine in the last section the automorphism groups of the Artin groups of rank 2. 
For such a group $A=A(\Gamma)$, we obtain that $\Out(A)$ is of order 2 if the label 
$m$ of the (unique) edge of $\Gamma$ is odd, and $\Out(A)$ is infinite and isomorphic 
to $(\Z \rtimes C_2) \times C_2$ if $m$ is even. 

\section{The semidirect product structures}\label{sect1}

We number the vertices of the Coxeter graphs of type $A_n$, $B_n$ and $D_n$ as shown in Figure \ref{fig1}.
The standard generators of $A(A_n)$ will be written $\al_1,\al_2,..,\al_n$, 
the standard generators of $A(B_n)$ will be written $\be_1,\be_2,..,\be_n$, 
and the standard generators of $A(D_n)$ will be written $\de_1,\de_2,..,\de_n$.
For the rest of the article we shall identify $A(A_{n-1})$ with the classical
$n$-string braid group $\Bn$ in the usual fashion. This will simplify notation throughout,
and be a constant remainder of the special role played by the braid group $\Bn$ in this story. 

\subsection{Definition of $\pi_B$, $\pi_D$, $s_B$, and $s_D$}

We wish to study the epimorphisms $\pi_B:A(B_n)\to \Bn$ and $\pi_D:A(D_n)\to \Bn$
defined by:

\[
\begin{aligned}
\pi_B(\be_1)=1\,,\hskip5mm \pi_B(\be_i)=\al_{i-1} \ \text{ for } i=2,3,..,n\,.\\
\pi_D(\de_1)=\pi_D(\de_2)=\al_1\,,\hskip5mm \pi_D(\de_i)=\al_{i-1} \ \text{ for } i=3,..,n\,.
\end{aligned}
\]

The epimorphism $\pi_B:A(B_n)\to\Bn$ admits a section $s_B:\Bn\to A(B_n)$ defined by:
\[
s_B(\al_i)=\be_{i+1}\ \text{ for } i=1,2,..,n-1\,.
\]
In particular, $A(B_n)$ may be written as a semidirect product $A(B_n)=\ker \pi_B \rtimes\Bn$.

In a similar fashion, the epimorphism $\pi_D:A(D_n)\to\Bn$ admits a section $s_D:\Bn\to A(D_n)$
defined by:

\[
s_D(\al_i)=\de_{i+1}\ \text{ for } i=1,2,..,n-1\,.
\]
and $A(D_n)$ may be written as a semidirect product $A(D_n)=\ker \pi_D \rtimes\Bn$.

\bigskip
We now describe how each of these semidirect products arise from a topologically defined
faithful action of the braid group $\Bn$ on a free group.
Both of these representations involved are classical. The first is quite famous and
due to Artin \cite{Ar1}.
The second arises from the monodromy action of the $n$-string braid group on the Milnor
fibre of the singularity of type $A_{n-1}$, and was shown to
be faithful by Perron and Vannier \cite{PV}. 

\begin{figure}[ht]
\begin{center}
\includegraphics[width=15cm]{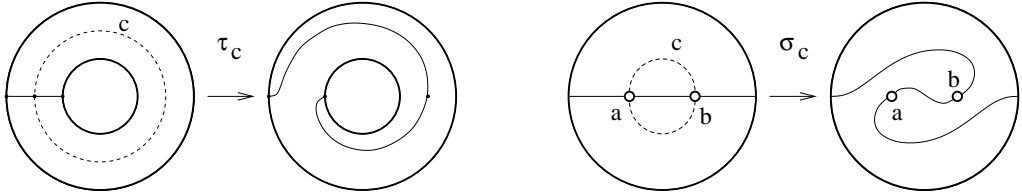}
\end{center}
\caption{Dehn twist $\tau_c$, and braid twist $\sig_c$.}\label{fig2}
\end{figure}

We recall the definition of a \emph{Dehn twist} homeomorphism $\tau_c:\Sigma\to \Sigma$ on
a simple closed curve $c$ in a surface $\Sigma$. Let $A$ denote an annular neighbourhood of $c$.
Then $\tau_c$ is any homeomorphism isotopic to one which is the identity on $\Sigma\setminus\text{int}(A)$ and
transforms the interior of $A$ as shown in Figure \ref{fig2}. 
If $c$ bounds a disk in $\Sigma$ and $a,b$ are points on $c$ which are
interchanged by $\tau_c$, then $\tau_c$ induces what we shall call a \emph{braid twist} homeomorphism $\sig_c$
of the punctured surface $\Sigma\setminus\{ a,b\}$. This homeomorphism exchanges
(neighbourhoods of) the two punctures, as shown in Figure \ref{fig2}.

\subsection{Artin's representation}\label{Artrep} 

\begin{figure}[ht]
\begin{center}
\includegraphics[width=8cm]{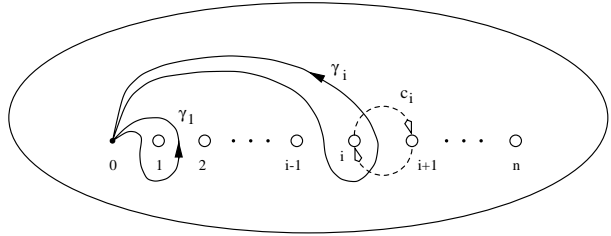}
\end{center}
\caption{The surface $\Sigma_B=\C\setminus\{ 1,2,..,n\}$, free group generators $u_i=[\ga_i]$, and braid twists
$\sig_i=\rho_B(\al_i)$ about $c_i$.}\label{fig3}
\end{figure}

Let $\Sigma_B$ denote the surface of Figure \ref{fig3},
namely $\C\setminus\{1,2,3,..,n\}$, and, for each $i=1,..,n-1$,
let $\sig_i$ denote the braid twist homeomorphism on the curve $c_i$ illustrated in Figure \ref{fig3}.
Take the origin $0\in\C$ as a basepoint and define loops $\ga_i$ at $0$ for $i=1,..,n$,
as shown in the figure.
Then $\pi_1(\Sigma_B,0)=F_n$ is freely generated by the
elements $u_i=[\ga_i]$. Artin's representation $\rho_B:\Bn\to\Aut(F_n)$ is defined by
$\rho_B(\al_i)=(\sig_i)_*$, for $i=1,2,..,n-1$.

\begin{prop}\label{Bsemi}
(1) Let $F_n$ denote the free group on $n$ generators $u_1,u_2,..,u_n$. Then Artin's representation
$\rho_B:\Bn\to\Aut(F_n)$ is well-defined and given algebraically by
\[
\rho_B(\al_i):
\begin{cases}
u_i\mapsto u_{i+1}\\
u_{i+1}\mapsto u_{i+1}^{-1}u_iu_{i+1}\\
u_j\mapsto u_j\hskip8mm j\notin\{ i,i+1\}\, .
\end{cases}
\] 

(2) We have $A(B_n)\cong F_n \rtimes_{\rho_B} \Bn$, where the projection onto the second factor is $\pi_B$
and the section $\Bn\hookrightarrow F_n\rtimes\Bn$ is just $s_B$. In particular, $\ker \pi_B$ is a free
group of rank $n$. 
\end{prop}

\Proof 
Part (1) is well-known and due to Artin \cite{Ar1}. It is easily enough checked that this
algebraic formulation gives a well-defined representation and corresponds to the topological description 
given above after identifying $F_n$ with $\pi_1\Sigma_B$ via $u_i=[\ga_i]$ for $i=1,2,..,n$.
We concern ourselves with part (2) of the proposition.

Let $\varphi:\{\be_1,..,\be_n\}\to F_n\rtimes\Bn$ be the function defined by $\varphi(\be_1)=u_1$
and $\varphi(\be_i)=\al_{i-1}$ for $i=2,..,n$. It is easily verified that $\varphi$ extends to a
homomorphism $\varphi:A(B_n)\to F_n\rtimes\Bn$. 

Let $\psi:\{u_1,..,u_n\}\cup\{\al_1,..,\al_{n-1}\}\to A(B_n)$ be defined by 
\[
\begin{aligned}
\psi(u_i)=\be_i\be_{i-1}...\be_2\be_1\be_2^{-1}...\be_{i-1}^{-1}\be_i^{-1}\hskip5mm \text{ for } i=1,..,n\,.\\
\psi(\al_i)=\be_{i+1}\hskip5mm \text{ for } i=1,..,n-1\,.
\end{aligned}
\]
Again, it is easily seen that $\psi$ induces a homomorphism $\psi:F_n\rtimes\Bn\to A(B_n)$.
(Hint: verify the relation $\varphi(\al_iu_iu_{i+1}\al_i^{-1})=\varphi(u_iu_{i+1})$, for all
$1\leq i<n$, by induction on $i$).  
Finally, one checks that $\varphi\circ\psi=\id$ and $\psi\circ\varphi=\id$.\endproof

\begin{prop}[Artin \cite{Ar1,Ar2}, Magnus \cite{Mag}]\label{ArtMag}
The representation $\rho_B:\Bn\to\Aut(F_n)$ is faithful and has image the set of automorphisms that
permute the conjugacy classes of the elements $u_1, u_2,..,u_n$ and fix the product $u_1u_2...u_n$. 
\noproof
\end{prop}


\subsection{Braid monodromy representation}\label{monodromy}

\begin{figure}[ht]
\begin{center}
\includegraphics[width=12cm]{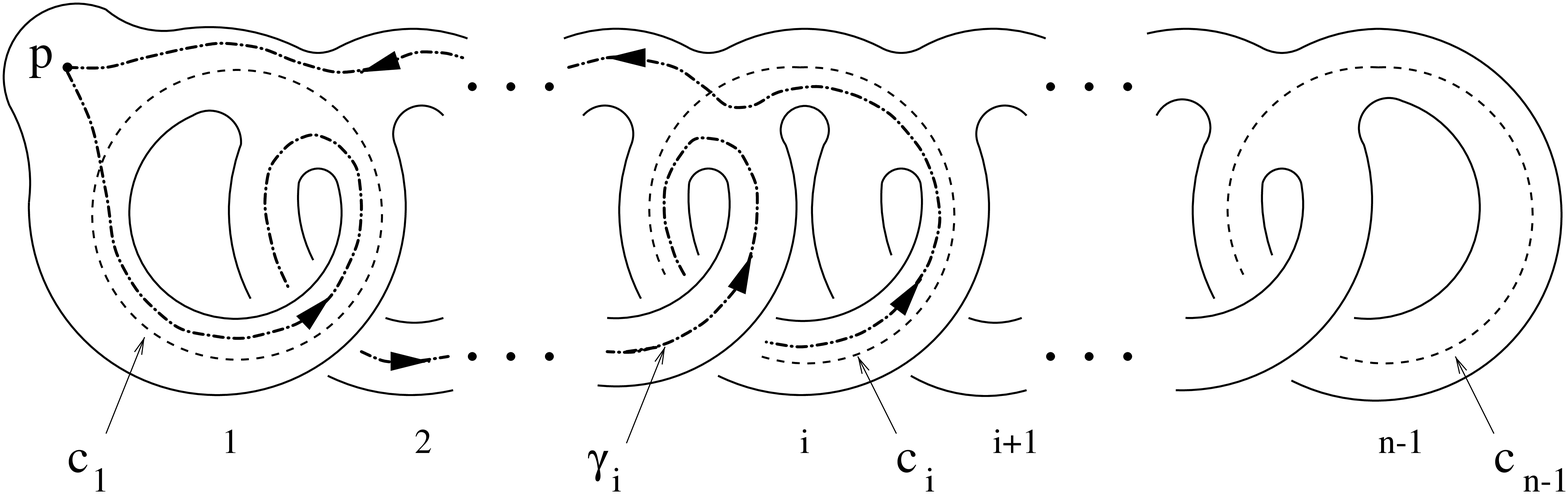}
\end{center}
\caption{The surface $\Sigma_D$, free group generators $v_i=[\ga_i]$, and Dehn twists
$\tau_i=\rho_D(\al_i)$ about $c_i$.}\label{fig4}
\end{figure}

Let $\Sigma_D$ denote the surface shown in Figure \ref{fig4} and, for $i=1,..,n-1$,
let $\tau_i$ denote the Dehn twist homeomorphism on the curve $c_i$ illustrated in the figure.
Choose a basepoint $p$ and loops $\ga_i$ at $p$ for $i=1,..,n-1$, as shown in the figure.
Then $\pi_1(\Sigma_D,p)=F_{n-1}$ is freely generated by the
elements $v_i=[\ga_i]$. The braid monodromy representation $\rho_D:\Bn\to\Aut(F_{n-1})$ is defined by
$\rho_D(\al_i)=(\tau_i)_*$, for $i=1,2,..,n-1$. 
  
\begin{figure}[ht]
\begin{center}
\includegraphics[width=12cm]{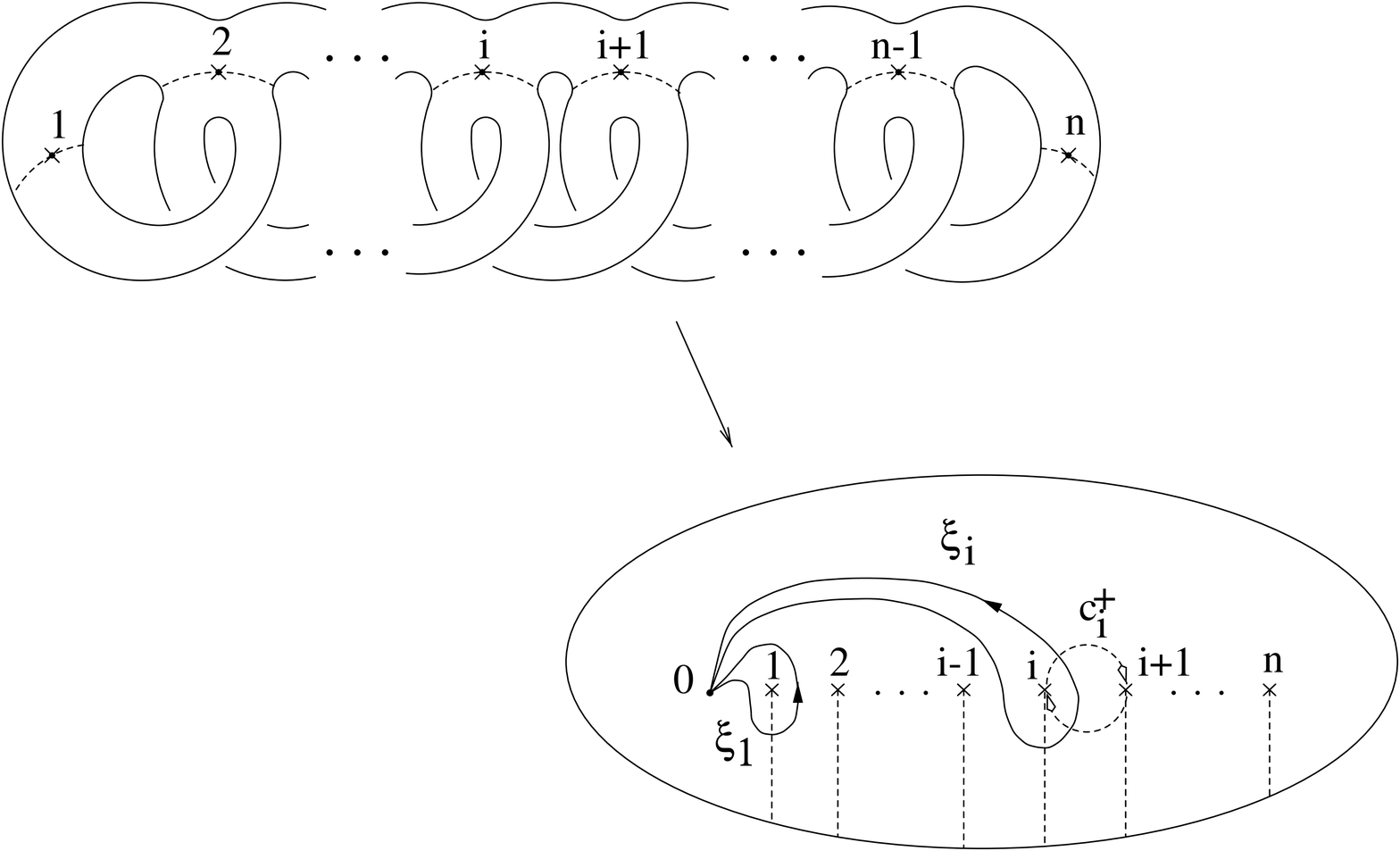}
\end{center}
\caption{The twofold branched cover $\Sigma_D\to \Sigma^+$, generators $x_i=[\xi_i]$ for $\star_n(C_2)$,
and braid twists $\sig^+_i=\rho^+(\al_i)$ about $c^+_i$.}\label{fig5}
\end{figure}

For what follows in the later sections, it is important to view the surface $\Sigma_D$ 
as a branched 2-fold cover of $\C$ branched over the set $\{1,2,..,n\}$, as shown in Figure \ref{fig5}. 
Note that the quotient orbifold $\Sigma^+$ has orbifold fundamental group $\star_n(C_2)$
(the free product of $n$ groups of order $2$).  
Let $\kappa: \star_n (C_2) \to C_2$ be the epimorphism which maps each free factor non trivially. 
Note that $\kappa$ is the unique epimorphism $\star_n (C_2) \to C_2$ which sends every element of 
order 2 to the generator of $C_2$, hence the kernel of $\kappa$ is a characteristic subgroup. 
Furthermore, $\ker \kappa =F_{n-1}$ is free of rank $n-1$ and the surface $\Sigma_D$ is the 
branched cover of $\Sigma^+$ corresponding to $F_{n-1} = \ker \kappa < \star_n (C_2)$.

We suppose that the basepoint $p$ in $\Sigma_D$ maps onto the origin $0\in \Sigma^+$.
The loops $\xi_i$ at $0$, for $i=1,..,n$, as shown in Figure \ref{fig5},
represent canonical generators $x_1,x_2,..,x_n$ for the group $\star_n(C_2)$.
The loops $\ga_i$ were chosen to be (homotopic to) lifts of the loops
$\xi_1\xi_{i+1}$, for $i=1,..,n-1$, and represent one choice of generators for
the free subgroup $F_{n-1}$.
A different choice will be used in a later section. The Dehn twist along the curve $c_i$
is a lift of the braid twist along the curve $c_i^+$ in $\Sigma^+$.
Thus $\rho_D$ extends to a representation $\rho^+:\Bn\to\Aut(\star_n C_2)$
which is induced from $\rho_B$ under the map $F_n\to\star_n C_2=F_n/(u_i^2, i=1,..,n)$. 
  
The following is simply a more explicit version of Perron and Vannier \cite{PV}, Corollaire 1, where it is 
observed that $A(D_n)$ is a semidirect product of a rank $n-1$ free group and the $n$-string
braid group.

\begin{prop}\label{Dsemi}
Suppose that $n\geq 4$.

(1) Let $F_{n-1}$ denote the free group on $n-1$ generators $v_1,v_2,..,v_{n-1}$. Then the braid monodromy
representation $\rho_D:\Bn\to\Aut(F_{n-1})$ is well-defined and given algebraically by:
\[
\rho_D(\al_1):
\begin{cases}
v_1\mapsto v_1\\
v_j\mapsto v_1^{-1}v_j\hskip8mm j\neq 1\,,
\end{cases}
\]
and, for $2\leq i\leq n-1$,
\[
\rho_D(\al_i):
\begin{cases}
v_{i-1}\mapsto v_i\\
v_i\mapsto v_iv_{i-1}^{-1}v_i\\
v_j\mapsto v_j\hskip8mm j\notin \{ i-1,i\}\, .
\end{cases}
\] 

(2) We have $A(D_n)\cong F_{n-1} \rtimes_{\rho_D} \Bn$, where the projection onto the second factor is $\pi_D$
and the section $\Bn\hookrightarrow F_{n-1}\rtimes\Bn$ is just $s_D$. In particular, $\ker\pi_D$ is a free
group of rank $n-1$. 
\end{prop}

\Proof 
One shows easily (by direct calculation) that this algebraic formulation of $\rho_D$ gives a well-defined
representation, and corresponds to the topological description above after identifying $F_{n-1}$ with $\pi_1\Sigma_D$
via $v_i=[\ga_i]$ for $i=1,..,n-1$.  
(The latter exercise is best carried out by looking at the braid twist action on $\pi_1(\Sigma^+)$.)
We turn to part (2).

Let $\varphi:\{\de_1,..,\de_n\}\to F_{n-1}\rtimes\Bn$ be the function defined by $\varphi(\de_1)=v_1\al_1$
and $\varphi(\de_i)=\al_{i-1}$ for $i=2,..,n$. One easily verifies that $\varphi$ extends to a
homomorphism $\varphi:A(D_n)\to F_{n-1}\rtimes\Bn$. 

Let $\psi:\{v_1,..,v_{n-1}\}\cup\{\al_1,..,\al_{n-1}\}\to A(D_n)$ be defined by 
\[
\begin{aligned}
\psi(v_i)=\de_{i+1}\de_i...\de_3\de_1\de_2^{-1}\de_3^{-1}...\de_i^{-1}\de_{i+1}^{-1}\hskip5mm
\text{ for } i=1,..,n-1\,.\\
\psi(\al_i)=\de_{i+1}\hskip5mm \text{ for } i=1,..,n-1\,.
\end{aligned}
\]
Again, it is easily seen that $\psi$ induces a homomorphism $\psi:F_{n-1}\rtimes\Bn\to A(D_n)$.
Finally, one checks that $\varphi\circ\psi=\id$ and $\psi\circ\varphi=\id$, and hence that $\varphi$
is an isomorphism.\endproof

\Remark The choice of the free basis for $F_{n-1}$ used here is the most convenient for the purposes of
this section. A different free basis, given by the elements $g_1:=x_1x_2=v_1$ and
$g_i:=x_ix_{i+1}=v_{i-1}^{-1}v_i$ for $i=2,3,..,n-1$, will be used later in Section \ref{sect3}.
With respect to this basis, $\rho_D$ is defined by:
\[
\rho_D(\al_i):
\begin{cases}
g_{i-1}\mapsto g_{i-1}g_i\\
g_{i+1}\mapsto g_i^{-1}g_{i+1}\\
g_j\mapsto g_j\hskip8mm j\notin\{ i-1,i+1\}\,.
\end{cases}
\]

\medskip

Faithfulness of the braid monodromy representation $\rho_D$
was established in \cite{PV} by essentially topological 
means (using the work of Birman and Hilden \cite{BH} on mapping class groups). We proceed now to
give a purely algebraic proof of faithfulness based on the above description. We first recall some
background and prove a lemma.

If $T$ is a subset of the standard generators of an Artin group $A=A(\Gamma)$, then the subgroup $A_T$ of
$A$ generated by $T$ is known as a \emph{standard parabolic subgroup} or \emph{special subgroup} of $A$.
By a well-known result of Van der Lek \cite{vdL}, each standard parabolic subgroup $A_T$ of $A(\Gamma)$ 
is canonically isomorphic to the Artin group $A(\Gamma_T)$ associated to the full subgraph $\Gamma_T$ 
of $\Gamma$ spanned by $T$.
 
\begin{lemma}\label{Lemma1}
Let $A$ be an Artin group of spherical type with standard generating set $S$. Let $T\subset S$,
$\al\in S\setminus T$ and $f,g\in A_T$. If $\al^{-1}f\al=g$, then $f=g\in A_{T\cap\al^\perp}$,
where $\al^\perp =\{\be\in S: m_{\al,\be}=2\}$. 
\end{lemma}

\Proof
In this proof we use the \emph{orthogonal normal forms} for elements
in a spherical type Artin group \cite{Ch1,Ch2}. Any spherical type Artin group admits 
a left-invariant lattice order $(A,<,\vee,\wedge)$ with positive cone the submonoid $A^+$ 
generated by the standard generators. Thus, for $x,y\in A$,
$x<y$ if and only if $x^{-1}y\in A^+$. This lattice order restricts nicely to
standard parabolic subgroups: $A^+\cap A_T$ is the submonoid generated by $T$ (written $A_T^+$)
for any $T\subset S$. As a result of this lattice structure, every element $x\in A$ has a unique
\emph{orthogonal normal form} $x=x_1^{-1}x_2$ where $x_1,x_2\in A^+$ and $x_1\wedge x_2=1$. In fact,
$x_1=(1\wedge x)^{-1}$ and $x_2=(1\wedge x)^{-1}x$. Clearly also,if $x\in A_T$ then both
$x_1,x_2\in A_T^+$.

Write $f=f_1^{-1}f_2$ and $g=g_1^{-1}g_2$ in orthogonal normal forms.
We have $f_1,f_2,g_1,g_2\in A_T^+$. The equality $\al^{-1}f\al=g$ implies that
$(f_1\al)^{-1}(f_2\al)=g_1^{-1}g_2$. Put $h=f_1\al\wedge f_2\al$. Since $g_1\wedge g_2=1$, we have 
$f_1\al=hg_1$ and $f_2\al=hg_2$. We now use the (nontrivial) fact that the submonoid $A^+$ of $A$
is presented as a monoid as follows:
\[
A^+ = \< S\ |\  w(\al,\be:m_{\al,\be})=w(\be,\al:m_{\al,\be})\text{ for all }\al,\be\in S\>^+\,.
\]
In other words, any two positive words representing the same element of $A^+$ are related
by a finite sequence of applications of the given relators. Using this, we see that, since $f_1\in A_T^+$
and $\al\in S\setminus T$, any positive word representing $f_1\al$ is of the form
$\ga_1\ga_2...\ga_p.\al.\de_1\de_2...\de_q$  
where $\ga_1,..,\ga_p\in T$ and $\de_1,..,\de_q\in T\cap\al^\perp$. Now, since $f_1\al=hg_1$ and $\al$ does
not appear in any positive word representing $g_1$, it follows that $g_1$ is
written $\de_1\de_2...\de_q$ with $\de_1,..,\de_q\in T\cap\al^\perp$.
Therefore $g_1\in A_{T\cap\al^\perp}^+$.
Similarly we have $g_2\in A_{T\cap\al^\perp}^+$ and therefore
$f=\al^{-1}g\al=g=g_1^{-1}g_2\in A_{T\cap\al^\perp}$. \endproof

The following is the ``$A_n$ case" of Theorem 1 of \cite{PV}. A completely analogous argument can also
be used to prove that Artin's representation $\rho_B:\Bn\to\Aut(F_n)$ is faithful.

\begin{prop}[Perron-Vannier \cite{PV}]\label{faithful}
The representation $\rho_D:\Bn\to\Aut(F_{n-1})$ is faithful.
\end{prop}

\Proof
We consider the spherical type Artin group $A=A(D_n)$. Recall the presentation
$A(D_n)=F_{n-1}\rtimes_{\rho_D}\Bn$ where $\Bn$ is identified with the standard parabolic subgroup
$A_{\{\de_2,..,\de_n\}}$ of $A$. To prove the proposition it suffices to show that,
if $f\in A_{\{\de_2,..,\de_n\}}$ commutes with $v_k$ for all $k=1,2,..,n-1$, then $f=1$.

Let $f\in A_{\{\de_2,\de_3,..,\de_n\}}$ such that $fv_1=v_1f$. Since $v_1=\de_1\de_2^{-1}$, we have 
\[\de_1^{-1}f\de_1=\de_2^{-1}f\de_2\in A_{\{\de_2,\de_3,..,\de_n\}}\,.\]
By Lemma \ref{Lemma1}, it follows that $f\in A_{\{\de_2,\de_4,..,\de_n\}}$.

We now show the following statement by induction on $k$ for $k\geq 3$:
\begin{description}
\item{($I_k$):} If $f$ commutes with $v_1,v_2,..,v_{k-1}$, 
then $f\in A_{\{\de_{k+2},\de_{k+3},..,\de_n\}}$.
\end{description}

Suppose $k=3$. Since $f$ commutes with $v_2$ and $v_2=\de_3v_1\de_3^{-1}$, we have that
$\de_3^{-1}f\de_3$ commutes with $v_1$. From the preceding argument we deduce that both $f$ and 
$\de_3^{-1}f\de_3\in A_{\{\de_2,\de_4,..,\de_n\}}$, and then, by Lemma \ref{Lemma1}, that
$f\in A_{\{\de_5,\de_6..,\de_n\}}$.

Suppose $k\geq 4$. Since $f$ commutes with $v_{k-1}$ and $v_{k-1}=\de_kv_{k-2}\de_k^{-1}$,
we have that $\de_k^{-1}f\de_k$ commutes with $v_{k-2}$. On the other hand, 
$\de_k^{-1}f\de_k$ also commutes with $v_1,v_2,..,v_{k-3}$ (since both
$\de_k$ and $f$ do).
By the induction hypothesis we therefore deduce that both
$f$ and $\de_k^{-1}f\de_k\in A_{\{\de_{k+1},\de_{k+2},..,\de_n\}}$, 
and then, by Lemma \ref{Lemma1}, that
$f\in A_{\{\de_{k+2},\de_{k+3},..,\de_n\}}$.

The assertion ($I_k$) with $k=n-1$ implies that if $f$ commutes with all $v_1,v_2,..,v_{n-2}$
then $f=1$. \endproof


\subsection{Comparison of the two representations}\label{Sect2.4}

Recall that, given a Coxeter graph $\Gamma$, the canonical map $A(\Ga)\to W(\Ga)$ is determined by 
adding the relations $\al^2=1$, for each standard generator $\al\in S$, to the standard presentation
of $A(\Gamma)$. Via this map, the semidirect product structures $A(B_n)=F_n\rtimes \Bn$ and
$A(D_n)=F_{n-1}\rtimes\Bn$ induce semidirect product structures on the associated Coxeter groups:
\[
W(B_n)=C_2^n\rtimes\Sn\, \text{ and }\ W(D_n)=C_2^{n-1}\rtimes\Sn\,
\]
respectively, where $\Sn$ denotes the symmetric group on $\{1, \dots,n\}$.
These product structures are, of course, well-known,
see \cite{Bou}. The group $W(B_n)$, sometimes called the \emph{signed permutation
group}, is simply the group of symmetries of an $n$-cube spanned by an orthonormal basis in $\R^n$:
the normal subgroup $C_2^n$ acts by change of signs of the $n$ coordinates, and $\Sn$ acts by permuting the
coordinates, and so acts on $C_2^n$ by permuting the direct factors. There exists a well-known embedding
$p: W(D_n)\to W(B_n)$ which commutes with the projections onto $\Sn$ and which sends $C_2^{n-1}$ onto the
kernel of the map $C_2^n\to C_2$ which is nontrivial on each factor.
It is tempting to wonder whether the inclusion $p$ is actually induced by an inclusion
$\phi:A(D_n)\to A(B_n)$ which commutes with the projections $\pi_B$ and $\pi_D$ onto $\Bn$.
This, however, is not the case.

\begin{prop}\label{noembed}
There is no embedding $\phi:A(D_n)\to A(B_n)$ such that $\pi_D=\pi_B\circ\phi$.
\end{prop}

The proof of Proposition \ref{noembed} makes use of the following lemma from Dyer and Grossman \cite{DG}. 
We first fix some notation: $F_n$ is the free group of rank $n$ generated by $u_1,u_2,..,u_n$ and $\Bn$ acts on 
$F_n$ via the representation $\rho_B$ as in Proposition \ref{Bsemi} -- for convenience we identify $\Bn$ with
its image under $\rho_B$, thus simplifying the notation.

For $x,y\in F_n$ we write $x\sim y$ if $x$ and $y$ are conjugate.

\begin{lemma}[Dyer-Grossman \cite{DG}, Lemma 18]\label{LemmaDyGr}
Let $x\in F_n$ such that $\al_i^2(x)\sim x$ for all $i=1,2,..,n-1$. Then $x$ is conjugate either
to a power of  $u_j$ for some $j\in\{1,2,..,n\}$ or to a power of $u_1u_2...u_n$.\noproof
\end{lemma}
 
\begin{corollary}\label{CorollaryCrucial}
Let $x \in F_n$ such that $\alpha_i(x) \sim x$ for all $i=1, \dots, n-1$. Then $x$ is conjugate to a power of $u_1u_2 \dots u_n$.
\end{corollary}

\Proof
By Lemma \ref{LemmaDyGr}, $x$ is conjugate either to a power of $u_j$ for some $j \in \{1, 2, \dots, n\}$ or 
to a power of $u_1u_2 \dots u_n$. If $x \sim u_j$, then we would have $\alpha_i(u_j) \sim u_j$ for 
all $i$, while in fact $\alpha_j(u_j)=u_{j+1}$ for $j<n$ and $\alpha_{j-1}(u_j) \sim u_{j-1}$ for $j>0$, 
a contradiction.
\endproof

\paragraph{Proof of Proposition \ref{noembed}.} 
We suppose that such an embedding $\phi:A(D_n)\hookrightarrow A(B_n)$ exists and shall henceforth identify
$A(D_n)$ with a subgroup of $A(B_n)$ via this map $\phi$. Recall that, by Proposition \ref{Bsemi},
$A(B_n)=F_n\rtimes \Bn$ where $F_n=\ker \pi_B$ is the free group on generators $u_1,u_2,..,u_n$,
$\Bn$ is the subgroup of $A(B_n)$ generated by $\be_2,\be_3,..,\be_n$,
and one has $\be_ix\be_i^{-1}=\rho_B(\al_{i-1})(x)=\al_{i-1}(x)$ for all $i=2,3,..,n$ and $x\in F_n$. 
 Recall also (Proposition \ref{Dsemi}) that
$A(D_n)=F_{n-1}\rtimes \Bn$ where $F_{n-1}=\ker \pi_D$ is the free group on generators $v_1,v_2,..,v_{n-1}$,
$\Bn$ is the subgroup of $A(D_n)$ generated by $\de_2,\de_3,..,\de_n$,
and one has $\de_iy\de_i^{-1}=\rho_D(\al_{i-1})(y)$ for all $i=2,3,..,n$ and $y\in F_n$. 
The fact that $\pi_D=\pi_B\circ\phi$ means that $F_{n-1}\subset F_n$ and, for all $i=2,3,..,n$, there is a $w_i\in F_n$
such that $\de_i=w_i\be_i$.

\paragraph{Step 1:} There exists $x_0\in F_{n-1}\setminus\{1\}$ such that $\al_i(x_0)\sim x_0$ for all
$i=1,2,..,n-1$.   

\[ 
\text{Define}\hskip15mm x_0 = 
\begin{cases}
v_1v_2^{-1}v_3v_4^{-1}...v_{n-2}^{-1}v_{n-1}\qquad &\text{if $n$ even,}\\
v_1v_2^{-1}v_3...v_{n-2}v_{n-1}^{-1}v_1^{-1}v_2v_3^{-1}...v_{n-2}^{-1}v_{n-1}\qquad &\text{if $n$ odd.}
\end{cases}\hskip15mm \,
\]
It is easily checked that $\de_ix_0\de_i^{-1}=\rho_D(\al_{i-1})(x_0)=x_0$. Therefore
\[
\al_{i-1}(x_0)=\be_ix_0\be_i^{-1}=w_i^{-1}\de_ix_0\de_i^{-1}w_i=w_i^{-1}x_0w_i
\hskip8mm\text{for all }i=2,..,n\,.
\]

For $x\in F_n$ we denote by $[x]$ the class of $x$ in $H_1(F_n)\cong\Z^n$.

\paragraph{Step 2:} We have $[x]=0$ for all $x\in F_{n-1}$.
\medskip

We show that $[v_i]=0$ by induction on $i$. Suppose that $i=1$. We have 
\[
 v_1^{-2}v_2 =\rho_D (\alpha_1^2) (v_2)=\de_2^2v_2\de_2^{-2}=(w_2\be_2)^2v_2(w_2\be_2)^{-2}=
w_2\al_1(w_2)\al_1^2(v_2)\al_1(w_2^{-1})w_2^{-1}
\]
and $\al_1^2$ acts trivially on $H_1(F_n)$, whence $[v_2]=-2[v_1]+[v_2]$, and so $[v_1]=0$.

Suppose now that $i\geq 2$ and $[v_{i-1}]=0$. Then
\[
\begin{aligned}
\,[v_i] = [\de_iv_{i-1}\de_i^{-1}]&=[\be_iv_{i-1}\be_i^{-1}]\hskip8mm
\text{since } \de_i=w_i\be_i \text{ with } w_i\in F_n\\
&=[\al_{i-1}(v_{i-1})]=\al_{i-1}([v_{i-1}])=0\,.
\end{aligned}
\]

\paragraph{End of proof.} Let $x_0$ be the element constructed in Step 1.
By Corollary \ref{CorollaryCrucial}, there exists some $k\in\Z$ such that $x_0\sim(u_1u_2..u_n)^k$. 
Step 2 now implies that 
\[
0=[x_0]=k([u_1]+[u_2]+..+[u_n])\,.
\]
Since the $[u_i]$ are independent nontrivial generators of $H_1(F_n)\cong\Z^n$, we must
have $k=0$, and therefore $x_0=1$, a contradiction.\endproof 

The fact that the pure braid group acts trivially on homology ($H_1(F_n)\cong\Z^n$) via Artin's
representation figures strongly in the work of Dyer and Grossman on the automorphisms of $\Bn$.    
We note here that the action of the braid group via the braid monodromy is somewhat more complicated
on the homology $H_1(F_{n-1})\cong \Z^{n-1}$. In fact, this action on $\Z^{n-1}$ turns out to be a
specialisation of the reduced Burau representation of $\Bn$.

\section{A topological interpretation for the semidirect products}\label{sect2}

Let $\Ga$ be a Coxeter graph of spherical type (i.e: such that the Coxeter group $W(\Ga)$
is finite).
Every finite Coxeter group has a canonical representation $\theta:W(\Ga)\hookrightarrow O(n,\R)$
where  $n$ is the number of vertices of $\Ga$, and where the standard generators are sent
to reflections.
It is well-known (since Brieskorn \cite{Br1}) that each Artin group $A(\Gamma)$ of spherical type
is isomorphic to the fundamental group of the space of regular orbits of the representation $\theta_\C$ of $W(\Ga)$
as a complex unitary reflection group, obtained simply by tensoring $\theta$ with $\C$. 
This regular orbit space may be easily described as the quotient of a complex hyperplane complement as follows.
 
Let $\Cal R$ denote the set of \emph{reflections} in $W(\Ga)$, that is the set of $r\in W(\Gamma)$ such that $\theta(r)$
is a reflection in a hyperplane, $H(r)$ say. (Note that every reflection is actually conjugate to a
standard generator of $W(\Ga)$).
The \emph{arrangement} of $\Ga$ is the set $\A(\Ga)=\{ H(r) : r\in\Cal R\}$ of reflecting hyperplanes of $W(\Ga)$.
The \emph{complexification} of $\A(\Ga)$ is the set $\A_\C(\Ga)=\{ H_\C=H\otimes\C : H\in\A(\Ga)\}$. Note that
$H_\C(r)=H(r)\otimes\C$ is simply the fixed set of $\theta_\C(r)$ in $\C^n$, and is a hyperplane in $\C^n$.
The \emph{complement} of $\A_\C(\Gamma)$ is the manifold
\[
M(\Ga)=\C^n\setminus \big( \bigcup_{H\in\A(\Ga)} H_\C\big)\ .
\]
The group $W(\Ga)$ acts freely on $M(\Ga)$, and the manifold $N(\Gamma)=M(\Ga)/W(\Ga)$ is the
regular orbit space referred to above. Brieskorn \cite{Br1} showed that $\pi_1N(\Ga)\cong A(\Gamma)$ and
Deligne later showed that  $N(\Gamma)$ is a $K(A(\Ga),1)$-space \cite{Del}.\footnote{In fact this was already known
to Brieskorn \cite{Br2} by more or less ad hoc means in all but 5 exceptional cases, including the cases we will 
consider in this article. Deligne's unified treatment was given in response to Brieskorn's
S\'eminaire Bourbaki and has significantly influenced the subsequent study of Artin groups.}
We refer the reader to Brieskorn's paper \cite{Br1} for an explicit description of the group isomorphism.

For types $B_n$, $D_n$ and $A_{n-1}$ we have the following (this information can be derived, for instance,
from the appendices of \cite{Bou}, where the associated root systems are laid out).
\[
\begin{aligned}
M(B_n) &= \{ (z_1,z_2,..,z_n)\in\C^n : z_i\neq\pm z_j\,\text{ for } 1\leq i\neq j\leq n\,, \\
&\qquad\text{ and } z_i \neq 0\,, \text{ for }1 \leq i\leq n\}\,,\hskip1cm\\ 
M(D_n) &= \{ (z_1,z_2,..,z_n)\in\C^n : z_i\neq\pm z_j\,,\text{ for } 1\leq i\neq j\leq n\}\,.\\
\end{aligned} 
\]
Let $D=\span\{(1,1,..,1)\}\subset\C^n$, and, for ${\bf z}\in \C^n$, let $[{\bf z}]$ denote 
the element of $\C^n/D$ represented by ${\bf z}$. Then
\[
M(A_{n-1}) = \{ [z_1,z_2,..,z_n]\in\C^n/D : z_i\neq z_j\,,\text{ for } 1\leq i\neq j\leq n\}\,.  
\]

The group $W(A_{n-1})=\Sn$ acts on $M(A_{n-1})$ by permutation of the coordinate axes,
$W(B_n)= (C_2)^n\rtimes\Sn$ acts on $M(B_n)$ by signed permutations of the coordinate axes
(i.e: $\Sn$ acts by permuting the coordinates and $(C_2)^n$ acts via
$(\varep_1,\varep_2,..,\varep_n).(z_1,z_2,..,z_n)=(\varep_1z_1,\varep_2 z_2,..,\varep_n z_n)$ for $\varep_i=\pm 1$),
and $W(D_n)$ acts on $M(D_n)$ also by signed permutations of the coordinates via the inclusion $W(D_n)<W(B_n)$
described in Subsection \ref{Sect2.4}.

We consider $M(B_n)\subset M(D_n)$ as submanifolds of $\C^n$
(i.e: $M(B_n)$ is simply the hyperplane complement $M(D_n)$ with additional hyperplanes removed,
namely the coordinate hyperplanes).
We also identify $W(D_n)$ as a subgroup of $W(B_n)$. The canonical actions of $W(B_n)$ and $W(D_n)$ on their respective 
hyperplane complements are therefore simultaneously induced by the action of $W(B_n)=(C_2)^n\rtimes\Sn$ on $\C^n$
by signed permutations of the coordinates.

\medskip

Let $\wtil P_D:M(D_n)\to M(A_{n-1})$ be the map defined by
$(z_1,z_2,..,z_n)\mapsto [z_1^2,z_2^2,..,z_n^2]$, and $\wtil P_B$ the restriction of this map
to the submanifold $M(B_n)\subset M(D_n)$. Note that these maps are equivariant with respect to the canonical projections
of $W(D_n)$ and $W(B_n)$ onto $\Sn=W(A_{n-1})$. Thus they induce, respectively, maps $P_D:N(D_n)\to N(A_{n-1})$ and 
$P_B:N(B_n)\to N(A_{n-1})$.

\begin{prop}\label{fibrations}
\begin{description}
\item[(B)]
{(i) The map $P_B:N(B_n)\to N(A_{n-1})$ is a locally trivial fibration. 
The fibre of $P_B$ over the orbit $\Sn.[\xi_1,\xi_2,..,\xi_n]$ is naturally homeomorphic to
$\C\setminus \{\xi_1,\xi_2,..,\xi_n\}$.\break
(ii)
The fibration admits a section $S_B$ and the subsequent
monodromy action of $\pi_1N(A_{n-1})$ on the fundamental group of the
fibre $\Sigma_B$ over the orbit $\Sn.[1,2,..,n]$ is precisely Artin's representation $\rho_B$.}\\
\item[(D)]
{(i) The map $P_D:N(D_n)\to N(A_{n-1})$ is a locally trivial fibration. 
The fibre of $P_D$ over any point is naturally homeomorphic to the twofold branched covering
of $\C$ branched over the set $\{\xi_1,\xi_2,..,\xi_n\}$.\break
(ii) The fibration admits a section $S_D$ and the subsequent
monodromy action of $\pi_1N(A_{n-1})$ on the fundamental group of the
fibre $\Sigma_D$ over the orbit $\Sn.[1,2,..,n]$ is precisely the braid monodromy representation $\rho_D$.}\\
\end{description}
\end{prop}

\Proof
We first treat the $B_n$ case. Let 
\[
\what M = M(B_n)/(C_2)^n = 
\{ (z_1,z_2,..,z_n)\in\C^n : z_i\neq z_j\,,\ z_i\neq 0\,,\text{ for } 1\leq i\neq j\leq n\}\,,
\]
where the covering $M(B_n)\to{\what M}$ is defined by the map $(z_1,z_2,..,z_n)\mapsto (z_1^2,z_2^2,..,z_n^2)$.
Now $N(B_n)={\what M}/\Sn$ and $P_B$ is induced by the map ${\what P}_B:{\what M}\to M(A_{n-1})$
given by ${\bf z}\mapsto [{\bf z}]$. Thus $\what P_B$ is simply the restriction
of the linear map $d:\C^n\to\C^n/D$ to the set $\what M$, and is easily seen to be a locally trivial fibration 
where the fibre over a point $[\bf \xi]$ is the set $({\bf\xi}+D)\setminus\{{\bf z}: z_i=0\text{ for some }i\}$,
thus naturally homeomorphic to $\C\setminus \{\xi_1,\xi_2,..,\xi_n\}$ by the map sending
$z\in\C\setminus\{\xi_1,..,\xi_n\}$ to the point $(\xi_1-z,\xi_2-z,..,\xi_n-z)$. 
Factoring out by the action
of $\Sn$ one obtains statement (i) for $P_B$. 

Furthermore, ${\what P}_B$ has a section, which can be described as follows.
For ${\bf z}=(z_1,z_2,..,z_n)\in\C^n$, define the  ``centre'' $c_{\bf z}:= \frac{1}{n}\sum z_i$
and ``radius'' $r_{\bf z}:=\text{max}\{|z_i-c_{\bf z}|\}$ of the set
$\{ z_1,z_2,..,z_n\}$. Then define ${\what S}_B:M(A_{n-1})\to\what M$ such that
${\what S}_B([{\bf z}])= (z_1+ b_{\bf z},z_2+b_{\bf z},..,z_n+b_{\bf z})$,
where $b_{\bf z}:=r_{\bf z}+1-c_{\bf z}$. 
Despite appearances, ${\what S}_B$ is well-defined and continuous.
Clearly, since ${\what S}_B$ is $\Sn$-equivariant,
it induces a section $S_B:N(A_{n-1})\to N(B_n)$ to the map $P_B$.
Note that the fibre over the orbit $\Sn.[1,2,..,n]\in N(B_n)$ is naturally homeomorphic to 
the surface $\Sigma_B$ of Figure \ref{fig3}, and that ${\what S}_B([1,2,..,n])=(1,2,..,n)$ and 
corresponds under
this homeomorphism to the basepoint at $0$ in $\Sigma_B=\C\setminus\{ 1,2,..,n\}$.

To see that the monodromy action is precisely Artin's braid action,
one simply recalls that the space of braids may be identified with
the space of loops in
$\{\,\{z_1,z_2,..,z_n\}\in\C^n/\Sn : z_i\neq  z_j\,,\text{ for } 1\leq i\neq j\leq n\}$ 
based at the point $\{ 1,2,..,n\}$. This permits an identification
of the braid group $\Cal B_n=A(A_{n-1})$ with $\pi_1N(A_{n-1})$
which happens to agree with that described in \cite{Br1}. 
(Note that the standard generator $\al_i$ of $A(A_{n-1})$ is identified with the
elementary braid which twists the $i^{\rm th}$ and $i+1^{\rm st}$ strands).

Now observe that the space $\what M^+ = M(D_n)/(C_2)^n$ is simply obtained from $\what M$ by restoring 
the coordinate hyperplanes. The map  $M(D_n)\to{\what M^+}$  (defined by $(z_1,z_2,..,z_n)\mapsto (z_1^2,z_2^2,..,z_n^2)$)
is a ramified covering with a singular set of degree 2, the union of restored hyperplanes. Thus $\what M^+$ should be thought
of as an orbifold. The fibration $\what P_B$ extends to an ``orbifold fibration" $\what P^+:\what M^+\to M(A_{n-1})$ with
the fibre over a point $[\bf \xi]$ naturally homeomorphic to the orbifold $\C$ with singular set 
$\{\xi_1,\xi_2,..,\xi_n\}$ of degree $2$ points (and therefore homeomorphic to the orbifold $\Sigma^+$ of
Figure \ref{fig5}).
This fibration descends to a fibration $P^+:N^+=\what M^+/{\Sn}\to N(A_{n-1})$
with the same fibre. Just as for $P_B$, the fibration $P^+$ admits a section $S^+$
and the monodromy action is simply that induced by Artin's braid action on the
punctured surface by replacing punctures by points of degree 2. This is precisely the
braid twist action of $\Sigma^+$ described in Section \ref{monodromy}.

It remains simply to observe that $N(D_n)$ is a twofold cover of $N^+$ (in fact $M(D_n)/\CE_{n-1}$ is clearly
a 2-fold cover of $\what M^+$) which fibres over $N(A_{n-1})$ as described in the statement (i), and moreover
to recall that the braid twist action on the orbifold $\Sigma^+$ lifts to the braid monodromy representation
on $\Sigma_D$. The existence of the section $S_D$ (a lift of $S^+$) is ensured simply by the fact that
the braid twist action on $\Sigma^+$ lifts.   
\endproof

\Remark The homotopy exact sequence implies that these fibrations give rise to semidirect product structures on the two
Artin groups $\pi_1N(B_n)\cong A(B_n)$ and $\pi_1N(D_n)\cong A(D_n)$.
Here we have not paid very much attention to explicit generators for the fundamental groups, however the careful reader
may verify that the canonical isomorphisms, described in Brieskorn \cite{Br1}, carry the
product structure coming from the fibration in each case precisely onto the semidirect product structures
given in Propositions \ref{Bsemi} and \ref{Dsemi}. 

\Remark In a sense, Proposition \ref{fibrations} is just a rephrasing of some previously known facts. 
As explained in  Daniel Allcock \cite{All},  the spaces $N(B_n)$, and $N^+$, may be regarded as the 
configuration spaces of $n$ unordered points in the orbifolds $\C\setminus\{0\}$, and $\C$ with a degree 2 singular point, respectively.
Thus $A(B_n)$ may be identified with the $n$-string braid group over a punctured plane (which is equally the subgroup
of the $(n+1)$-string braid group in which the first string is pure), and $A(D_n)$ with an index
two subgroup of the braid group over a plane with a degree 2 point. The fibrations observed above correspond simply
to deleting the orbifold features (the puncture or the degree 2 point) and Artin's representation appears quite naturally
from this point of view as well. We note that this braid picture for type $B_n$ had been previously observed \cite{Lam,Cr}
and, as Allcock pointed out, was already implicit in \cite{Br2}. In fact the fibration of type $D_n$ in Proposition
\ref{fibrations} is also implicit in \cite{Br2} where Brieskorn observes that the map $f:\C^n\to\C^{n-1}$ defined by
$(z_1,z_2,..,z_n)\mapsto (z_1^2-z_n^2,z_2^2-z_n^2,..,z_{n-1}^2z_n^2)$ restricts to a
locally trivial (differentiable) fibration of $M(D_n)$. This fibration is simply the map $\wtil P_D$, where
the image $f(M(D_n))$ is just $M(A_{n-1})$ expressed in the inhomogeneous coordinates obtained by setting the last
coordinate to $0$. We should point out that Allcock also gave similar ``braid picture'' interpretations of the
infinite families of affine type Artin groups. The picture for $A(\wtil A_n)$ has
already been used in \cite{CP}, and there is certainly potential
for these to be explored more fully.   


\section{Automorphisms preserving the fibre}\label{sect3}

This section is largely inspired by the paper of Dyer and Grossman \cite{DG} which gave the first
proof, using essentially algebraic techniques, that $\Out(\Bn)$ is the group of two elements
(the nontrivial outer automorphism being that which 
changes the sign of each standard braid group generator). While we are not yet able to replicate 
that result for the groups $A(B_n)$ and $A(D_n)$, we are able to prove that there are very few
automorphisms of these groups which leave invariant the kernel of the map $\pi_B$ and $\pi_D$
respectively (that is, the normal subgroups in the semidirect product structures described in the
previous sections). The question remains as to whether these kernels are characteristic subgroups.
On the other hand, we are able to show that $\Out( (\star_n C_2) \rtimes \Bn)$ is also of order 2.

In our investigation in this section of the automorphisms of the group $A(B_n)$ we shall always 
assume $n \ge 3$. The case $n=2$, which will be treated in Section \ref{chap5} together with the 
other Artin groups of rank 2, is manifestly different. For example, the group $A(B_2)$ has a 
non-inner automorphism derived from the non-trivial automorphism of the graph $B_2$, and such an 
automorphism does not exist if $n \ge 3$. In fact, the group $\Out (A(B_2))$ is infinite (see 
Theorem \ref{thm51}) while we expect that the group $\Out(A(B_n))$ is finite if $n \ge 3$. 


\subsection{Automorphisms of $A(B_n)$ leaving $F_n$ invariant}

Recall the presentation $A(B_n)=F_n\rtimes{\rho_B}\Bn$ of Section \ref{sect1}. 
For notational convenience, we shall identify elements of $\Bn<A(B_n)$
with the corresponding elements of $\Inn(A(B_n))$ via their action by conjugation 
(and using the fact that $\rho_B$ is faithful).
We call these elements \emph{braid automorphisms}.
As in Section \ref{sect1}, let $F_n$ be freely generated by $u_1,u_2,..,u_n$ so that $\rho_B$ is defined
as in Proposition \ref{Bsemi}. For simplicity we write $u_0:=u_1u_2u_3...u_n$, and note that
all braid automorphisms leave $u_0$ fixed.

\begin{figure}[ht]
\begin{center}
\includegraphics[width=15cm]{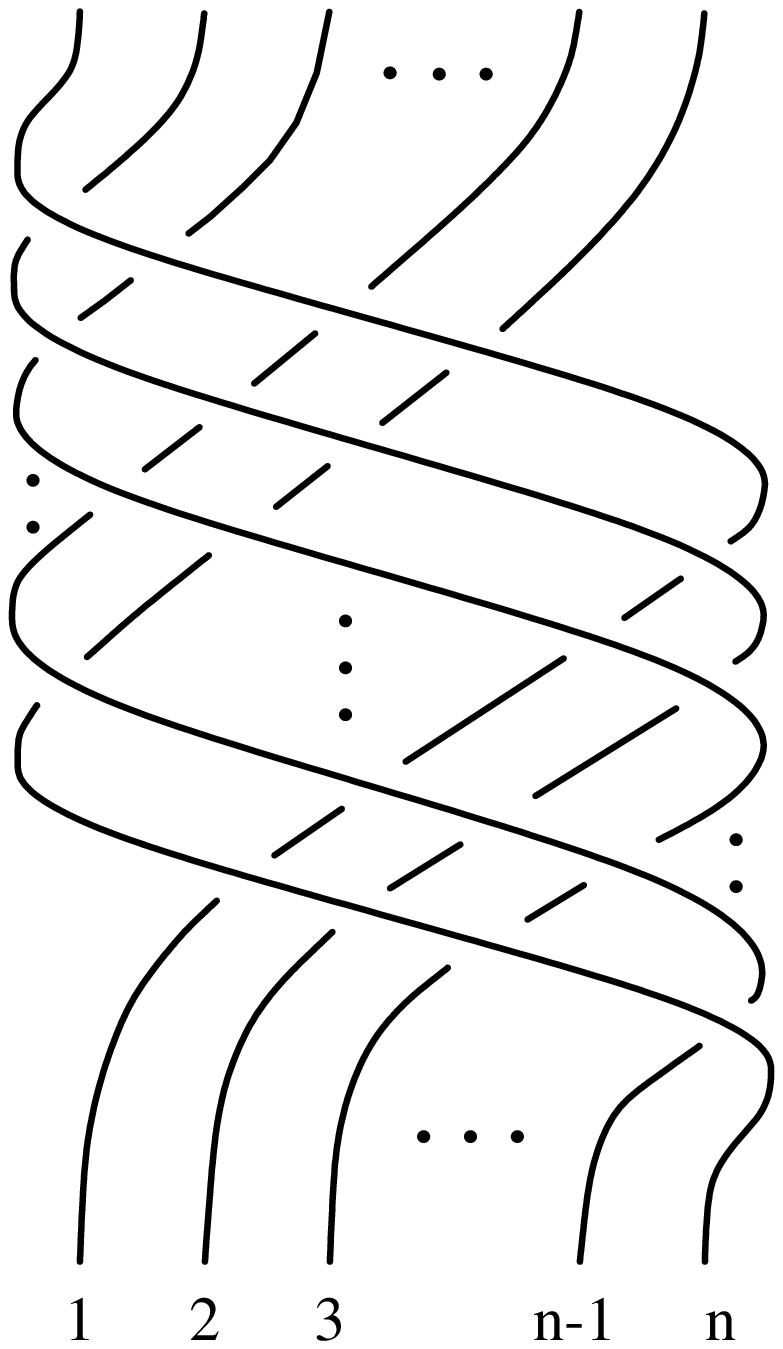}
\end{center}
\caption{The braid $\zeta=(\al_1\al_2\al_3\cdots\al_{n-1})^n$ in $\Bn$.}\label{fig6}
\end{figure}

Let $\zeta\in\Bn$ denote the braid of Figure \ref{fig6}.
\footnote{Note that braids are drawn from top to bottom, so that their action (on $\pi_1\Sigma_B$ 
or $\pi_1\Sigma^+$)
should be visualised by moving loops drawn in the bottom level surface 
up to the top level surface via a continuous path of loops.}
It is known that $\zeta$ generates the centre
of $\Bn$ (it is the square of Garside's so-called fundamental element), see \cite{Gar}. 
Also, it is easily seen that $\zeta$ acts on $A(B_n)$ by conjugation by the element $u_0^{-1}$,
and it is a straightforward exercise to check that the centre of $A(B_n)$ is generated by the
element $u_0\zeta$.

Note that there is an obvious bijective correspondence between automorphisms of $A(B_n)$ which
fix the subgroup $\Bn$ while leaving the subgroup $F_n$ invariant, and  automorphisms of $F_n$
which are equivariant with respect to the $\Bn$ action.
The proof of the following result is essentially contained in the work of Dyer and Grossman \cite{DG} --
 see the  proof of their Theorem 19.

\begin{prop}\label{Bequivautos}
Let $n\geq 3$, and let $\Bn$ act on $F_n$  
via Artin's representation $\rho_B$.
Then the group of $\Bn$-equivariant automorphisms of $F_n$ is generated by $\zeta$.
\end{prop}

\Proof
We first note that the group $F_n^\Bn$ of elements of $F_n$ left fixed by $\Bn$ is generated
by $u_0$: if $x$ is a nontrivial element of $F_n^\Bn$ then, by Corollary \ref{CorollaryCrucial},
$x$ is conjugate
to a nontrivial power of $u_0$, say $x=wu_0^rw^{-1}$.
But then the commutator $[w,u_0^r]$ lies in $F_n^\Bn$
and so is conjugate to a power of $u_0$. However, since $u_0$ maps to an infinite order element
under abelianisation of $F_n$ we must have $[w,u_0^r]=1$, in which case $w$ must be a power of $u_0$ and
$x\in\<u_0\>$.

Let $\varphi$ be a $\Bn$-equivariant automorphism of $F_n$. It follows from the above statement
that $\varphi$ must leave invariant the cyclic group $\< u_0\>$. Thus $\varphi(u_0)=u_0^\nu$ for some
$\nu\in\{\pm 1\}$.

We now use the fact that, for all $i=1,2,..,n-1$ and
all $j=0,1,2,..,n$, we have $\al_i^2(u_j)\sim u_j$, and therefore by $\Bn$-equivariance 
$\al_i^2(\varphi(u_j))=\varphi(\al_i^2(u_j))\sim\varphi(u_j)$. 
It follows from Lemma \ref{LemmaDyGr}, that
$\varphi$ must permute the $2n$ conjugacy classes of the elements
$u_1^{\pm 1},..,u_n^{\pm 1}$ (while respecting the pairs $\{ u_i, u_i^{-1}\}$).
 Note that the elements
$u_0^{\pm 1},u_1^{\pm 1},..,u_n^{\pm 1}$ represent distinct conjugacy classes in $F_n$, and
none of these elements are proper powers. The conjugacy classes of $u_0$ and $u_0^{-1}$ are already 
either fixed or interchanged.

So, for $j=1,2,..,n$, we have 
$\varphi(u_j)\sim u_{\sig(j)}^{\nu_j}$ for $\nu_j\in\{\pm 1\}$,
where $\sig$ denotes a permutation in $\Sn$. 
Moreover, after abelianising $F_n$ to $\Z^n$, the equation
$\varphi(u_0)=\varphi(u_1)\varphi(u_2)...\varphi(u_n)$  implies that $\nu_j=\nu$ for every $j=1,2,..,n$.
It now follows, by Proposition \ref{ArtMag}, that there exists a braid $\be\in\Bn$ such that 
$\be^{-1}\varphi(u_j)=u_j^\nu$ for $0\leq j\leq n$.

Now, if $\nu=-1$, we would have
\[
u_1^{-1}u_2^{-1}\cdots u_n^{-1}=\be^{-1}\varphi(u_1u_2\cdots u_n)
=\be^{-1}\varphi(u_0)=u_0^{-1}=u_n^{-1}\cdots u_2^{-1}u_1^{-1}\,,
\]
which is impossible. Therefore $\varphi$ agrees on $F_n$ with the braid automorphism $\be$.
But then, since $\varphi$ is $\Bn$-equivariant and since the monodromy representation is faithful,
$\be$ must lie in the centre of $\Bn$ which is generated by $\zeta$.
\endproof

\begin{thm}\label{Bautos}
Let $n\geq 3$ and consider the presentation of $A(B_n)=F_n\rtimes_{\rho_B}\Bn$ as a semidirect product.
The group $\Aut(A(B_n),F_n)$ of automorphisms of $A(B_n)$ which leave invariant the subgroup $F_n$
is generated by the inner automorphisms of $A(B_n)$ and the automorphism $\ep_n$ which simply inverts
each of the standard generators. Thus
\[
\Aut(A(B_n),F_n)=\Inn(A(B_n))\rtimes\< \ep_n \>\cong\big( A(B_n)/Z(A(B_n))\big)\rtimes C_2\,,
\]
where $Z(A(B_n))$ denotes the centre of $A(B_n)$. In particular, $\Out(A(B_n),F_n) \cong C_2$.
\end{thm}

\Proof 
Let $\varphi\in \Aut(A(B_n),F_n)$, and let $\ov{\varphi}$ denote the automorphism induced on 
the quotient $\Bn=A(B_n)/F_n$. 
Using the theorem of Dyer and Grossman \cite{DG}, we may suppose,
up to a multiplication of $\varphi$ by a braid automorphism, and by $\ep_n$ if necessary,
that $\overline{\varphi}$ is the trivial automorphism of $\Bn$. 
In other words, there exists some function $k:\Bn\to F_n$, written $\ga\mapsto k_\ga$,
such that
\[
\varphi(\ga)= k_\ga\,\ga\,,\hskip5mm \text{ for all } \ga\in\Bn\,.
\] 
In particular, for all $\ga\in\Bn$ and all $x\in F_n$, we have 
\[
\varphi\circ\ga(x)=\varphi(\ga)\varphi(x)\varphi(\ga)^{-1}=k_\ga\,\ga\,\varphi(x)\,\ga^{-1}k_\ga^{-1}
\sim \ga\circ\varphi(x)\,,
\]
where $\sim$ denotes conjugacy in $F_n$.
In other words, the action of $\varphi$ on the conjugacy classes of $F_n$ is $\Bn$-equivariant. 

Note that, for every $\ga\in\Bn$, we have $\ga(\varphi(u_0))\sim\varphi(\ga(u_0))=\varphi(u_0)$.
By Corollary \ref{CorollaryCrucial} and the fact that
$u_0$ is not a proper power, we then have $\varphi(u_0)\sim u_0^\nu$, for $\nu\in\{\pm 1\}$.
But then, modifying $\varphi$ by an inner automorphism coming from $F_n$,
we may suppose that $\varphi(u_0)=u_0^\nu$. (Note that, under this modification, $\ov{\varphi}$
remains trivial, however the function $k$ may change).

As a consequence, we now claim that $\varphi$ actually fixes every element of $\Bn$ (and so its 
 action on $F_n$ is genuinely $\Bn$-equivariant). All braid automorphisms leave $u_0$ fixed.
So, for any $\ga\in\Bn$ we have
\[ 
u_0^\nu=\varphi(u_0)=\varphi(\ga u_0 \ga^{-1})=k_\ga\,\ga\,u_0^\nu\,\ga^{-1}k_\ga^{-1}
=k_\ga\,u_0^\nu k_\ga^{-1}\,.
\]
That is, $k_\ga$ commutes with $u_0$ and so must be a power of $u_0$, for 
every $\ga\in\Bn$. In fact, $k:\Bn\to\< u_0\>\cong \Z$ must be a homomorphism since  
$\varphi(\ga\be)=k_\ga\,\ga\, k_\be\,\be=k_\ga k_\be\,\ga\be$,
using once again the fact 
that braid automorphisms fix powers of $u_0$.
Moreover, in view of the braid relations $\al_i\al_{i+1}\al_i=\al_{i+1}\al_i\al_{i+1}$,
the only homomorphisms $\Bn\to\Z$ are multiples of the length homomorphism $\ell:\Bn\to\Z$ which
sends each standard generator to 1. In other words, there exists an $m\in\Z$ such that
$\varphi(\ga)=u_0^{m\ell(\ga)}\ga$, for all $\gamma \in \Bn$. Finally, we observe that 
the centre of $A(B_n)$, which is generated by the element $u_0\zeta$, is left invariant by
any automorphism of $A(B_n)$. So $\varphi(u_0\zeta)=(u_0\zeta)^{\pm 1}$. 
But since $\varphi(u_0\zeta)=u_0^{\nu}\, u_0^{m\ell(\zeta)}\,\zeta$ we must have 
$m\ell(\zeta)=1-\nu$. But, since $\ell(\zeta)=n(n-1)\geq 6$, this is only possible if $m=0$,
in which case $\varphi$ fixes every element of $\Bn$ as claimed.

It now follows from Proposition \ref{Bequivautos} that $\varphi$ is a power of the
central braid automorphism $\zeta$, and hence an inner automorphism of $A(B_n)$.
\endproof


\subsection{Automorphisms of $A(D_n)$ leaving $F_{n-1}$ invariant}

We proceed now to derive the analogous result for the type $D_n$ Artin group.
The proof of Theorem \ref{Bautos} is largely transportable to this case. Thus, the main part of
our work will be in establishing the following analogue of Proposition \ref{Bequivautos}.

\begin{thm}\label{autos}
Let $n\geq 4$, and identify $\Bn$ with a subgroup of $Aut(F_{n-1})$ via the
braid monodromy representation $\rho_D$.
Then the group of $\Bn$-equivariant automorphisms of $F_{n-1}$ is generated by $\zeta$.
\end{thm}

Let $K=\< x_1,x_2,..,x_n\mid x_i^2=1, i=1,2,..,n\>$
denote the group $\pi_1\Sigma^+=\star_n(C_2)$, where
the generators $x_i$ are as described in Subsection \ref{monodromy}.
The braid group $\Bn$ acts on $K$ (via braid twists of the degree 2 orbifold points) in such a way
that, for $i=1,2,..,n-1$,
\[
\al_i : \begin{cases}
x_i\mapsto x_{i+1}\\
x_{i+1}\mapsto x_{i+1}x_ix_{i+1}\\
x_j\mapsto x_j\text{ if } j\not\in\{i,i+1\}\\
\end{cases}
\]  

The braid monodromy action (via $\rho_D$) is just the restriction of this action to the characteristic
index $2$ free subgroup $F_{n-1}<K$, namely the kernel of the map of $K$ onto $C_2$ which
maps each generator $x_i$ nontrivially. In this section we shall use the following set of free generators for
$F_{n-1}$: $\{g_j=x_jx_{j+1} : j=1,2,..,n-1\}$. 
 
We first study the action of $\Bn$ on the whole of $K$.
For $u,v\in K$ we write $u\sim v$ to mean that
$u$ and $v$ are conjugate in $K$, and write $[u]$ for the conjugacy class of $u$ in $K$.
Note that any element in $K$ is uniquely represented by a
\emph{reduced form}, a word in $x_1,..,x_n$ with all exponents $+1$. In addition, every conjugacy class in
$K$ has a representative whose reduced form is \emph{cyclically reduced},
i.e: every cyclic shift of the word is also reduced (equivalently, the word is reduced and 
does not begin and end with the same letter). Such a representative is unique up to cyclic shifts of its
reduced form, and shall be called a \emph{shortest} element
of its conjugacy class. The same observations hold in the case of reduced forms with
respect to free products of arbitrary groups.
The following three lemmas shall be also used in the next subsection, so we shall assume $n\ge 2$
in their statement in place of $n\ge 4$.
 
\begin{lemma}\label{first}
Let $1\leq j\leq n-1$, and let $K^{\<\al_j\>}$ denote the subgroup of $K$ consisting
of those elements left fixed by $\al_j$.
\begin{description}
\item{(i)} $K^{\<\al_j\>}$ is generated by the elements $\{ g_j=x_jx_{j+1}$, $x_1,..,x_{j-1}$,
$x_{j+2},..,x_n\}$.
\item{(ii)} If $w\in K$ is such that $\al_j(w)\sim w$, then there is a shortest element
of $[w]$ lying in $K^{\<\al_j\>}$.
\end{description} 
\end{lemma}

\Proof 
(i) The group $K$ may be written $C\star D$ where $C=\<x_j,x_{j+1}\>$ and $D=\<x_i:i\neq j,j+1\>$.
Clearly $\al_j$ fixes every element of $D$ and leaves $C$ invariant.
So $K^{\<\al_j\>} = C^{\<\al_j\>}\star D$.

For simplicity we write $\al=\al_j$, $x=x_j$ and $y=x_{j+1}$. 
We have $\al(x)=y$ and $\al(y)=yxy$, thus $\al(xy)=xy$ and
$\al(yx)=yx$. Any element of $C=\<x\>\star \<y\> \cong C_2 \star C_2$ may be written in normal form 
with respect to the free product and has one of the following forms:
$(xy)^k, (yx)^k=(xy)^{-k}, (xy)^kx$, or $y(xy)^k$, for some $k\in\N$.
Now $\al((xy)^kx)=(xy)^ky\neq (xy)^kx$, and $\al(y(xy)^k)=yxy(xy)^k\neq y(xy)^k$. Thus $C^{\<\al\>}$
is generated by $xy=g_j$, and $K^{\<\al_j\>}=\<g_j\>\star D$ is as claimed.

(ii) 
We may clearly assume $w\neq 1$
and may choose $w$ to be any shortest element of its conjugacy class.
Thus, either $w\in C$, or $w\in D$, or $w$ may be chosen to have reduced form 
$W=c_1d_1c_2d_2..c_md_m$ with respect to the free product $C\star D$, where the $c_i$
are nontrivial elements of $C$ and the $d_i$ are nontrivial elements of $D$, and $m\geq 1$.

Suppose first that $w\in C$. 
Since it is shortest in its conjugacy class, $w$ must be represented
by a cyclically reduced word in the letters $x,y$. That is, $w$ is either
$(xy)^k$ or $(yx)^k=(xy)^{-k}$, for some $k\in\N$. But then $w$ lies in $K^{\<\al_j\>}$.

If $w\in D$ then, a fortiori, $w$ is  an element of $K^{\<\al_j\>}$.

Finally we suppose that $w$ has reduced form $W=c_1d_1c_2d_2..c_md_m$ (with respect to $C\star D$), 
in which case it follows that $W':=\al(c_1)d_1...\al(c_m)d_m$ is the reduced form for $\al(w)$
with respect to $C\star D$ . Since $\al(w)\sim w$, and since both $W$ and $W'$ are
cyclically reduced, we must have that $W'$ is obtained from $W$ by
a cyclic shift. That is, there exists an $r\in\{1,2,..,m\}$ such that
$\al(c_i)=c_{i+r}$ and $d_i=d_{i+r}$ for all $i=1,2,..,m$ (where indices are taken mod $m$).
But then $\al^{2m}(c_i)=c_i$ for all $i=1,2,..,m$.
Note that $\al^2$ acts on both $x$ and $y$ by conjugating by the element $yx$. But then
$\al^{2m}(c_i)=c_i$ only if $c_i$ is a power of $xy$. It follows that $w$ must lie in $K^{\<\al_j\>}$.
\endproof

The following lemma should be compared with Lemma 18 of \cite{DG} (see Lemma \ref{LemmaDyGr}) and, more
particularly, with Corollary \ref{CorollaryCrucial}. 

\begin{lemma}\label{deltakey}
\begin{description}
\item {(i)} The group $K^\Bn$ of elements of $K$ left fixed by every $\ga\in\Bn$ is the cyclic group
generated by the element $\delta:=x_1x_2...x_n$. 
\item{(ii)} If $w\in K$ is such that $\al_j(w)\sim w$ for all $j=1,2,..,n-1$, then
$w$ is conjugate to an element of $K^\Bn$ (so conjugate to a power of $\de$).
\end{description}
\end{lemma}

\Proof
The case $n=2$ is contained in Lemma \ref{first}. We assume $n \ge 3$.

(i) Let $w \in K^\Bn$ and let $u_1 u_2 ..u_k$ be the reduced form for $w$ (each $u_i \in \{x_1, x_2, 
.., x_n\}$).

(1) Suppose that $u_i=x_j$ where $1<j<n$. Then, by Lemma \ref{first}, the inclusions $w \in K^{\< 
\alpha_{j-1} \>}$ and $w \in K^{\< \alpha_j \>}$ imply that $1<i<k$ and
$u_{i-1} u_i u_{i+1}=x_{j-1} x_j x_{j+1}$, or $x_{j+1} x_j x_{j-1}$.

(2) Suppose that $u_i=x_1$. Then, by Lemma \ref{first}, the inclusion $w \in K^{\< \alpha_1 \>}$ 
implies that either $u_{i-1}=x_2$ or $u_{i+1}=x_2$. Moreover, we cannot have both identities, 
$u_{i-1}=x_2$ and $u_{i+1}=x_2$, otherwise, by (1), $x_3x_2x_1x_2x_3$ would be a subexpression of 
$u_1 u_2 ..u_k$ (recall that $n \ge 3$), and then the inclusion $w \in K^{\< \alpha_1 \>}$ would 
contradict Lemma \ref{first}.

(3) Suppose $u_i=x_n$. Then, by Lemma \ref{first}, the inclusion $w \in K^{\< \alpha_{n-1} \>}$ 
implies that either $u_{i-1}=x_{n-1}$ or $u_{i+1}=x_{n-1}$. Moreover, we cannot have both 
identities, $u_{i-1}=x_{n-1}$ and $u_{i+1}=x_{n-1}$, otherwise, by (1), $x_{n-2}x_{n-1}x_nx_{n-
1}x_{n-2}$ would be a subexpression of $u_1 u_2 ..u_k$, and then the inclusion $w\in K^{\< 
\alpha_{n-1} \>}$ would contradict Lemma \ref{first}.

Clearly, (1), (2), and (3) imply that $w$ is a power of $\delta$. On the other hand, 
it is easily verified that $\delta$ is indeed fixed by every braid.
\medskip

(ii) Let $w \in K$ such that $\alpha_j(w) \sim w$ for all $j=1, ..., n-1$, and let $u_1 u_2 
..u_k$ be a cyclically reduced word representing an element of $[w]$ (each $u_i \in \{x_1, x_2, 
.., x_n\}$).

(1) Suppose $u_i=x_j$, where $1<j<n$. Then, by Lemma \ref{first}, we have $u_{i-1} u_i u_{i+1}= 
x_{j-1} x_j x_{j+1}$, or $x_{j+1} x_j x_{j-1}$. (This time, the indices of the $u_i$'s are
considered modulo $k$. So $u_0=u_k$ and $u_{k+1}=u_1$.)

(2) Suppose $u_i=x_1$. Then, by Lemma \ref{first}, either $u_{i-1}=x_2$ or $u_{i+1}=x_2$, and we 
cannot have both identities, $u_{i-1}=x_2$ and $u_{i+1}=x_2$.

(3) Suppose $u_i=x_n$. Then, by Lemma \ref{first}, either $u_{i-1}=x_{n-1}$ or $u_{i+1}=x_{n-
1}$, and we cannot have both identities, $u_{i-1}=x_{n-1}$ and $u_{i+1}=x_{n-1}$.

Now, (1), (2), and (3) clearly imply that $u_1 u_2 ..u_k$ is a power of $\delta$ up to a cyclic 
shift.   
\endproof

We will only need to use part (i) of the following lemma, but we state part (ii) anyway for the sake of
completeness.

\begin{lemma}\label{second} 
Let $S$ denote the set of standard generators of $\Bn$.
Choose a sequence $1\leq i_1<i_2<\cdots <i_k <n$ and let $T\subset S$ be the set 
$T=S\setminus\{\al_j : j=i_1,i_2,..,i_k\}$. Let $K^{\<T\>}$ denote the subgroup of $K$ consisting
of those elements left fixed by every element of $T$. Then:
\begin{description}
\item {(i)} $K^{\< T\>}=\<\de(1,i_1),\de(i_1+1,i_2), \de(i_2+1,i_3),..,\de(i_k+1,n)\>$, where, for $i<j$,
$\de(i,j):=x_ix_{i+1}\ldots x_{j-1}x_j$.
\item{(ii)} If $w\in K$ is such that $\al_j(w)\sim w$ for all $\al_j\in T$, then
$w$ is conjugate to an element of $K^{\< T\>}$.
\end{description}
\end{lemma}

\Proof
(i) Write $i_0=1$ and $i_{k+1}=n$. Consider the decomposition $K=K_0 \star K_1 \star \dots \star 
K_k$, where $K_r$ is the subgroup of $K$ generated by $x_{i_k+1}, \dots, x_{i_{k+1}}$. Let $w 
\in K^{\<T\>}$, and let $w=w_1 w_2 \dots w_p$ be the reduced form of $w$ with respect to this 
decomposition, namely, $w_i \in K_{\mu(i)}$ for some $\mu(i) \in \{0, 1, \dots, k\}$, and 
$\mu(i+1) \neq \mu(i)$ for all $i=1, \dots, p-1$. Now, by exactly the argument of Lemma 
\ref{deltakey} (i), we deduce that, if $w_i \in K_r$ (namely, $r=\mu(i)$), then $w_i$ is a power 
of $\delta (i_r+1,i_{r+1})$.

(ii) This follows by a similar extension of Lemma \ref{deltakey} (ii).
\endproof

We will also make use of the following two lemmas in the proof of Theorem \ref{autos}.

\begin{lemma}\label{third}
Let $\varphi$ be an automorphism of the free group $F\<x,y\>$ such that 
$\varphi(y)=y^\nu$ for some $\nu\in\{\pm1\}$. Then
$\varphi(x)=y^kx^\varep y^l$ for some $k,l\in\Z$ and $\varep\in\{\pm 1\}$.  
\end{lemma}

\Proof
Write $\varphi(x) = y^kuy^l$ where $k,l\in\Z$ and $u$ is a nontrivial reduced word whose
first and last letters are either $x$ or $x^{-1}$. We will use the fact that, since $\varphi$ is
an automorphism, the group $F\< x,y\>$ is generated by $y$ and $u$. 
The reduced word $u$ is written uniquely in the form $u_0u_1u_0^{-1}$ where the subword $u_1$ is cyclically reduced
and nontrivial. It follows that, for any $r\in\Z\setminus\{ 0\}$, $u^r$ has reduced form $u_0u_1^ru_0^{-1}$.
Thus the reduced form for $u^r$ begins and ends in $x$ or $x^{-1}$ and has length at least $\length(u)$.

Let $w(s,t)$ be a reduced word in the letters $s,t$ involving at least one $t$. The above observation
shows that the length of $w(y,u)$ is greater than or equal to the length of $u$. So, $x$ can never lie 
in the subgroup generated by $y$ and $u$ unless $u$ is of length $1$. But then 
$\varphi(x)$ must be of the form $\varphi(x) = y^k x^\varep y^l$ with $\varep \in \{ \pm 1 \}$.
\endproof

\begin{lemma}\label{fourth}
Let $w(s,t)$ denote a freely reduced word in the letters $s,t$. Then, in the free group $F_2\< x,y\>$,
we have the relation $w(x,xy).w(y,xy)=1$ if and only if $w(s,t)$ is the trivial word.
\end{lemma}

\Proof
Clearly the relation holds if $w(s,t)$ is the trivial word. Suppose then that the relation holds, and 
let $u(s,t)$ be the freely reduced word in $s,t$ such that $u(x,y)=w(x,xy)$.
Then we have
\[
u(x,y)y^{-1}u(y,x)y=w(x,xy)y^{-1}w(y,yx)y=w(x,xy)w(y,xy)=1\,.
\]
But then $u(y,x)y=yu(x,y)^{-1}$ and, comparing lengths, it is clear that either $u(x,y)=1$,
or the word $u(y,x)$ is obtained from the word $u(x,y)^{-1}$ by a cyclic shift involving a single
letter (namely a $y$ or $y^{-1}$). But the latter is impossible. So $w(x,xy)=u(x,y)=1$. But since
mapping $s\mapsto x$ and $t\mapsto xy$ defines an isomorphism $F_2\<s,t\>\to F_2\<x,y\>$, $w(x,xy)=1$
only if $w(s,t)$ is the trivial word.
\endproof 

\paragraph{Proof of Theorem \ref{autos}.}
The proof of Theorem \ref{autos} falls into two cases, depending on whether $n$ is even or odd.
Observe that the braid $\zeta$ acts on $F_{n-1}$ by conjugation by the element $\delta^{-1}$ where
$\delta :=x_1x_2..x_n\in K$,
and that $\delta\in F_{n-1}$ if and only if $n$ is even.

\paragraph{Proof of Theorem \ref{autos}, case $n$ even.} 

Let $\varphi:F_{n-1}\to F_{n-1}$ be some
$\Bn$-equivariant automorphism. Our objective is to prove that $\varphi$ is a power of $\zeta$.
The first observation is that $\varphi$ (\emph{and} its inverse) must leave invariant
the fixed subgroup $F_{n-1}^\Ga=K^\Ga\cap F_{n-1}$ for any $\Ga <\Bn$.
That is, $\varphi$ restricts to an automorphism of $F_{n-1}^\Ga$.

 By Lemma \ref{deltakey}, $F_{n-1}^\Bn=K^\Bn=\<\de\>$. Therefore,
by the above observation, $\varphi$ restricts to an automorphism of $\<\de\>$,
hence $\varphi(\de)=\de$ or $\de^{-1}$.
We also have, from Lemma \ref{second}, that the two elements $g_1$ and $\de$ freely generate the subgroup
$F_{n-1}^{\< T\>}=K^{\< T\>}$, where $T=\{\al_1\}\cup\{\al_3,..,\al_{n-1}\}$.
Thus $\varphi$ restricts to an automorphism
$\varphi:F\<g_1,\de\>\to F\<g_1,\de\>$. Lemma \ref{third} now applies to show that
$\varphi(g_1)=\de^kg_1^\varep\de^l$ for some $k,l\in \Z$ and $\varep\in\{\pm 1\}$.

We now observe that, for each $i=2,3,..,n-1$, $g_i=\ga_i(g_1)$ for some braid $\ga_i\in\Bn$. Thus by equivariance
of $\varphi$ and the fact that $\de$ is fixed by all braids we have
\[
\varphi(g_i)=\de^kg_i^\varep\de^l \hskip5mm \text{ for fixed } k,l\in\Z\text{ and } \varep=\pm 1\,,
\text{ and for all $i$.}
\]
Now, since $n$ is even, we have $\de=g_1g_3..g_{n-1}\in F_{n-1}$. Therefore,
\begin{equation}\label{deleqn}
\varphi(\de)=\de^kg_1^\varep\de^{l+k}g_3^\varep\de^{l+k}\ldots\de^{l+k}g_{n-1}^\varep\de^l=\de^\nu \hskip5mm
\text{ for some }\nu=\pm 1\,.
\end{equation}
Abelianising $F_{n-1}$ to $\Z^{n-1}$ the relation (\ref{deleqn}) 
becomes $(m(l+k)+\varep-\nu)[\de]=0$, where $n=2m\geq 4$.
There are two cases: 
\begin{description}
\item {(i)} $\varep=\nu$ and therefore $l+k=0$, or
\item{(ii)} $\varep=-\nu$ and therefore $m.(l+k)= 2\nu$, so that, necessarily, $n=4$ and $l+k=\nu$.
\end{description} 

In case (i), the possibility $\varep=\nu=-1$ leads to a contradiction, for it implies that 
$g_1g_3...g_{n-1}=g_{n-1}...g_3g_1$. But then $\varep=1$, and $\varphi$ evidently just acts on $F_{n-1}$ by
conjugation by $\de^k$. That is, $\varphi=\zeta^{-k}$ as required.

\begin{figure}[ht]
\begin{center}
\includegraphics[width=15cm]{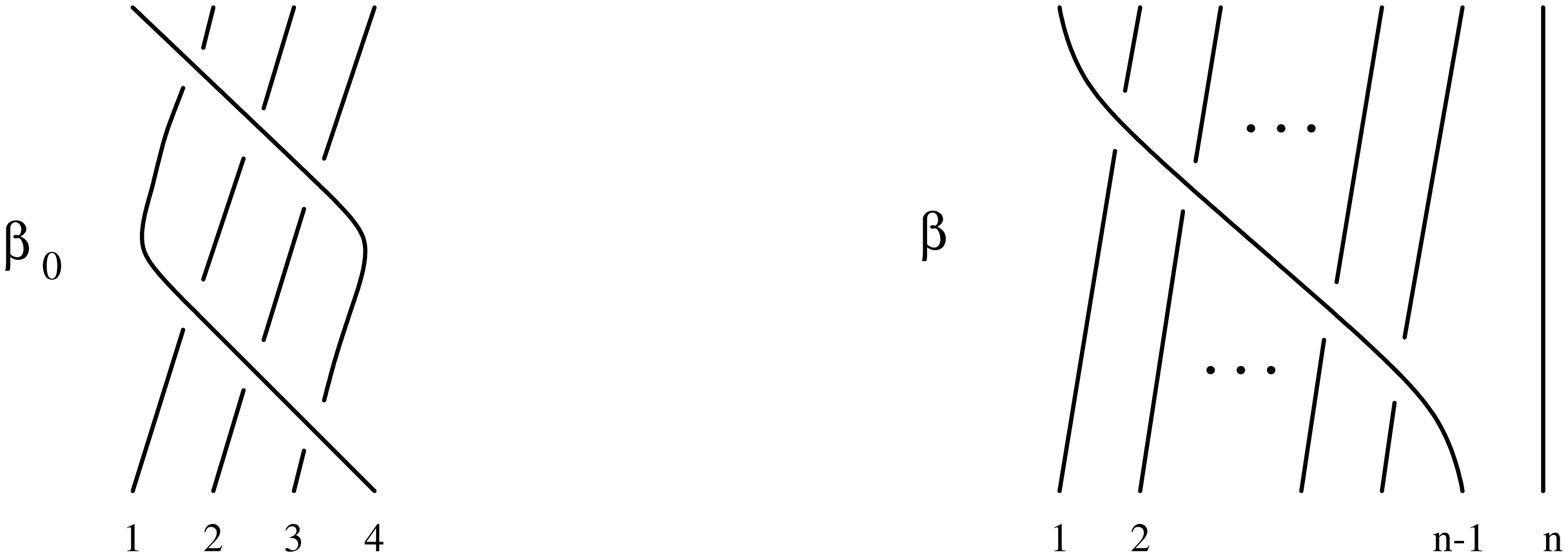}
\end{center}
\caption{The braids $\be_0=(\al_1\al_2\al_3)^2$ and $\be=\al_1\al_2...\al_{n-2}$.}\label{fig7}
\end{figure}

In case (ii), $n=4$, $\de=g_1g_3$, $l+k=\nu=-\varep$, and (\ref{deleqn}) simplifies to
$g_1^\varep(g_1g_3)^{-\varep} g_3^\varep=1$, which is possible if and only if $\varep=-1$.
But then, $\zeta^{k}\varphi$ is the automorphism $F_3\to F_3$ sending
$g_i$ to $g_i^{-1}\delta$ for $i=1,2,3$. That is $g_1\mapsto g_3$, $g_2\mapsto g_2^{-1}g_1g_3$
and $g_3\mapsto g_3^{-1}g_1g_3$. Let $\be_0\in\Cal B_4$ denote the braid $(\al_1\al_2\al_3)^2$,
shown in Figure \ref{fig7}.
Then one can easily check that $\zeta^{k}\varphi$ is realised by the action of $\be_0$,
which is a contradiction since $\be_0$ does not lie in the centre of $\Cal B_4$
and so is not $\Cal B_4$-equivariant (while $\zeta^{k}\varphi$ is).
\endproof
  
\paragraph{Proof of Theorem \ref{autos}, case $n$ odd.}

As in the previous case, let $\varphi:F_{n-1}\to F_{n-1}$ be some
$\Bn$-equivariant automorphism. Again, we shall prove that $\varphi$ is a power of $\zeta$,
using the observation  that $\varphi$ restricts to an automorphism of $F_{n-1}^\Ga=K^\Ga\cap F_{n-1}$
for any $\Ga <\Bn$. 

Define $x:=x_1x_2x_3\ldots x_{n-1}$. Then $x=g_1g_3\ldots g_{n-2}\in F_{n-1}$ and $\de=xx_n$.
Note that this time $\de$ does not lie in $F_{n-1}$, however the element
$z:=\de^2=xx_nxx_n$ does. Therefore $F_{n-1}^\Bn=\< z\>$
(since, by Lemma \ref{deltakey}, $K^\Bn=\<\de\>$).
Since $\varphi$ restricts to an automorphism of this group, we have $\varphi(z)=z$ or $z^{-1}$.

\paragraph{Claim 1.} {\it After multiplication of $\varphi$ by a power of $\zeta$ we may suppose
that $\varphi(x)=x$.}

\Proof
We consider the subset of braid generators $T=\{\al_1,\al_2,\ldots ,\al_{n-2}\}$.
By Lemma \ref{second}, $K^{\<T\>}=\<x\>\star \<x_n\> \cong \Z \star C_2$, and $F_{n-1}^{\<T\>}$ is therefore
freely generated by the two elements $x$ and $\what x:=x_nxx_n$. Note that $x\what x=xx_nxx_n=\de^2=z$,
so that $F_{n-1}^{\<T\>}$ is also freely generated by $x$ and $z$. Thus $\varphi$ restricts to an
automorphism of $\<x,z\> \cong F_2$. Since $\varphi(z)=z$ or $z^{-1}$, Lemma \ref{third} applies to give 
\[
\varphi :\begin{cases}
x\mapsto z^kx^\varep z^l\\
z\mapsto z^\nu
\end{cases}
\]
for some $k,l\in\Z$ and $\varep,\nu\in\{\pm 1\}$.


Notice that $\zeta(x)=\de^{-1}x\de =\what x$ (since $\de=xx_n$).
Consequently, by equivariance,
we have that $\varphi(\what x)=z^k\what x^\varep z^l$. Thus
\begin{equation}\label{zeqn}
\varphi(z)=\varphi(x\what x)=z^kx^\varep z^{l+k}\what x^\varep z^l=z^\nu\,. 
\end{equation}
Abelianising $F_2$ to $\Z^2$, this gives $(2(k+l)+\varep-\nu)[z]=0$. There are two cases: 
\begin{description}
\item {(i)} $\varep=\nu$ and $l+k=0$, or
\item{(ii)} $\varep=-\nu$ and $(l+k)=-\varep$.
\end{description} 

In case (i) above, (\ref{zeqn}) becomes $x^\varep\what x^\varep=(x\what x)^\varep$ which is
only possible if $\varep=1$ in which case $\varphi(x)=z^kxz^{-k}=\zeta^{-2k}(x)$.

In case (ii) above, (\ref{zeqn}) becomes $x^\varep(x\what x)^{-\varep}\what x^\varep=1$ which
is only possible if $\varep=-1$, and so $l+k=1$. In this case
$\varphi(x)=z^kx^{-1}z^l=\zeta^{-2k}(x^{-1}z)=\zeta^{-2k}(\what x)=\zeta^{-2k+1}(x)$, since
$\zeta(x)=\what x$.
\endproof

We now restrict our attention to the subgroup
$H=\< g_1,g_2,g_3,...,g_{n-2}\>\star \< x_n\> \cong F_{n-2} \star C_2$ of $K$.
Let $\what\cdot:H\to H$ denote the involution which is conjugation by $x_n$. (This is consistent with
the definition of the element $\what x$ already introduced). Then we see that
$H\cap F_{n-1}= A \star B$ where $A$ denotes the group $\< g_1,g_2,..,g_{n-2}\>\cong F_{n-2}$ and 
$B:= \what A$. Thus $B$ is freely generated by the elements
$\what g_1,\what g_2,\what g_3,..,\what g_{n-2}$. (Note: $A\star B$ is just the kernel of the map 
$H\to C_2$ induced by the quotient $K\to K/F_{n-1}$). 
Given a free product of groups $G=G_1\star G_2$, each element $g\in G$ has a unique expression in  
\emph{reduced form} with respect to the free product decomposition, namely an expression
\[
g=u_1u_2u_3\cdots u_k
\]
where each syllable $u_i$ is a nontrivial element of either $G_1$ or $G_2$ and where two consecutive
syllables do not belong to the same factor $G_1$ or $G_2$. The number $k$ shall be referred to as the
\emph{syllable length} of $g$ and written $||g||$.

\medskip

Consider the set of braid generators $T=\{\al_1\}\cup\{\al_3,..,\al_{n-2}\}$. By
Lemma \ref{second} 
\[
K^{\< T\>}=\<g_1,g_3g_5..g_{n-2}\>\star \<x_n\> = \<g_1,x\>\star \<x_n\> \cong F_2 \star C_2\,.
\]
Let $A_1:=\<g_1,x\>$ and $B_1:=\what A_1 =\<\what g_1,\what x\>$. Then
$K^{\< T\>}\cap F_{n-1}= A_1\star B_1$, and $\varphi$ restricts to an automorphism of
this group. In particular, $\varphi(g_1)\in A_1\star B_1$. Let 
\[
W=w_1w_2w_3\cdots w_k
\]
be the reduced form for $\varphi(g_1)$ with respect to $A_1\star B_1$.
Thus the syllables $w_i$ are nontrivial elements coming alternately from the groups $A_1$ and $B_1$.
Since $A_1<A$ and $B_1<B$, this is equally a reduced form with respect to $A\star B$. Also $W$ is
nontrivial since $\varphi(g_1)$ is nontrivial.

Now let $\be\in\Bn$ denote the braid $\al_1\al_2\cdots\al_{n-2}$ shown in Figure \ref{fig7}. 
Observe that $\be(g_i)=g_{i+1}$ if $i\leq n-3$ and $\be(x)=x$ (so that
$\be(g_{n-2})= (g_2g_4\cdots g_{n-3})^{-1}g_1g_3\cdots g_{n-2}$ using the fact
that $x=g_1g_3\cdots g_{n-2}$). Thus $\be$ leaves the subgroup $A$ invariant. Moreover, since
$\be(x_n)=x_n$, the braid $\be$ leaves the whole of $A\star B$ invariant respecting the free product
structure. Thus, for instance, $\be^r(W)=\be^r(w_1)\be^r(w_2)\cdots\be^r(w_k)$ is a reduced form
with respect to $A\star B$ for the element $\be^r(\varphi(g_1))=\varphi(g_{r+1})$, for all $r\in \Z$.
Moreover, these reduced forms all have similar structure:  a syllable $\be^r(w_i)$
comes from the factor $A$ if and only if $w_i$ does.

\paragraph{Claim 2.}
{\it If $\varphi(x)=x$, then $\varphi(g_1)\in A_1$.}

\Proof
We continue with the notation introduced in the preceding paragraphs.
By the hypothesis  $\varphi(x)=x$, we may use the fact that $\varphi(x)$ 
has syllable length 1 with respect to $A\star B$.
However, since $g_j=\be^{j-1}(g_1)$, for $j=1,2,..,n-2$ (this is NOT true for $j=n-1$), and
since $\varphi$ acts $\Bn$-equivariantly, we have
\[
\varphi(x)=\varphi(g_1g_3g_5\cdots g_{n-2})= W\be^2(W)\be^4(W)\cdots \be^{n-3}(W)\,.
\]

If $||W||$ is even then $W\be^2(W)\be^4(W)\cdots \be^{n-3}(W)$ is already in reduced form 
with respect to $A\star B$, and has syllable length at least $4$ ($n\geq 5$ and $||W||\geq 2$)
which contradicts the uniqueness of reduced forms in a free product.

We may assume therefore that $||W||$ is odd, and let $m=\frac{||W||+1}{2}$. Then $W$ is of the
form $U^{-1}w_mV$ where $U$ and $V$ are reduced forms of the same syllable length
with respect to $A\star B$ and whose first syllables
are of the same \emph{type} (i.e: they come from the same factor $A$ or $B$). 

Write $V=v_1v_2\cdots v_{m-1}$ and $U=u_1u_2\cdots u_{m-1}$ as reduced forms. 
If $V\be^2(U)^{-1}\neq 1$ then there is a last $i\in\{1,2,..,m-1\}$ such that
$v_i\neq\be^2(u_i)$, in which case $V\be^2(U)^{-1}$ has reduced form
\[
M= v_1v_2\cdots v_{i-1}.\,\omega.\,\be^2(u_{i-1}^{-1})\cdots \be^2(u_2^{-1})\be^2(u_1^{-1})
\]
where $\omega= v_i\be^2(u_i)^{-1}$ is a nontrivial element of the same free factor as $v_i$.
This implies that 
\[
U^{-1}w_mM\be^2(w_m)\be^2(M)\cdots\be^{n-5}(M)\be^{n-3}(w_m)\be^{n-3}(V)
\]
is a reduced form for $\varphi(x)$ of syllable length at least $3$ ($n\geq 5$).
But this again contradicts uniqueness of the reduced form.  
Thus we may suppose that $V=\be^2(U)$.

Both $V$ and $U$ represent elements of $A_1\star B_1$. Thus $V=\be^2(U)$ represents
an element of $(A_1\star B_1) \cap \be^2(A_1\star B_1)$ which may also be written 
$\<g_1,\what g_1,x,\what x\>\cap \<g_3,\what g_3,x,\what x\>$. 
 Except in the case $n=5$ where $x=g_1g_3$,
the elements $g_1,g_3,x,\what g_1,\what g_3,\what x$ form a free system.
Therefore, in the case $n\geq 7$, the intersection of these two rank 4 free groups is just
$\<x,\what x\>\cong F_2$ and is therefore fixed elementwise by $\be$. Since $V$ lies in this intersection,
we now have $U=\be^{-2}(V)=V$. Thus $\varphi(x)$
has reduced form $U^{-1}\omega U$ with middle syllable $\omega :=w_m\be^2(w_m)\cdots\be^{n-3}(w_m)$.
Note that $\omega\neq 1$ since it is conjugate to $\varphi(x)\neq 1$.
But since $||\varphi(x)||=1$ we must have $U=1$. Therefore $m=1$ and $\varphi(g_1)=Uw_mU^{-1}=w_1$ and
clearly lies in $A_1$ since otherwise we would have $\omega\in B$ contradicting the fact that
$\omega=\varphi(x)=x\in A$. This completes the proof of the Claim in the case $n\geq 7$.

In the case $n=5$, we have $\varphi(x)= U^{-1}w_m\be^2(w_m)\be^2(V)$ (since $V=\be^2(U)$).  
Since it is a syllable, $w_m$ is a nontrivial element of either $A_1$ or $B_1$.
Suppose that $w_m\in A_1\setminus\{1\}$. 
Then, there is a nontrivial freely reduced word $w(s,t)$ in the letters $s,t$, such that 
$w_m\be^2(w_m)= w(g_1,x)w(g_3,x)$, where in this case $x=g_1g_3$. 
But then, by Lemma \ref{fourth}, $\omega := w_m\be^2(w_m)\neq 1$
and  $U^{-1}\omega\be^2(V)$ is therefore a reduced form for $\varphi(x)$. The fact that
$||\varphi(x) ||=1$ then tells us that $U=V=1$, so that $\varphi(g_1)=w_m\in A_1$, as required.
If, however, $w_m\in B_1$, then the same argument shows that $\varphi(g_1)=w_m\in B_1$ and
$\varphi(x)=w_m\be^2(w_m)\in B$, which would contradict $\varphi(x)=x$.
\endproof

By Claim 1, we may suppose that $\varphi(x)=x$.
In that case, by Claim 2, $\varphi(A_1)<A_1$.
But, since $\varphi^{-1}(x)=x$, applying Claim 2 to $\varphi^{-1}$ also gives
$\varphi^{-1}(A_1)< A_1$, and so $\varphi$ restricts to an automorphism of $A_1$.
Lemma \ref{third} now applies to give $\varphi(g_1)=x^kg_1^\varep x^l$ for some $k,l\in\Z$ 
and $\varep=\pm 1$. By equivariance of $\varphi$ with respect to $\be$, we   
in fact have $\varphi(g_j)=x^kg_j^\varep x^l$ for all $j=1,2,..,n-2$ (but NOT necessarily for $j=n-1$).
But then 
\[
x=\varphi(x)=x^kg_1^\varep x^{l+k}g_3^\varep x^{l+k}\ldots x^{l+k}g_{n-2}^\varep x^l
\]
Abelianising $F_{n-1}$ to $\Z^{n-1}$ yields the equation $\frac{(n-1)}{2}(l+k)+\varep=1$.
We either have 
\begin{description}
\item {(i)} $\varep=1$ and $l+k=0$, in which case $\varphi(g_i)=x^kg_ix^{-k}$ for $i\leq n-2$, or
\item{(ii)} $n=5$, $l+k=1$ and $\varep=-1$, in which case $\varphi(g_i)=x^kg_i^{-1}x^{1-k}$ for $i=1,2,3$.
\end{description} 

Now let $y=g_2g_4\cdots g_{n-1}$. By a similar analysis
(labelling the $x_i$ in the reverse order.) we arrive at the conclusion that
there exist $m,l\in\Z$ such that either
\begin{description}
\item {(iii)} $\zeta^m\varphi(g_i)=y^lg_iy^{-l}$  for all $i=2,3,..,n-1$, or
\item {(iv)} $\zeta^m\varphi(g_i)=y^lg_i^{-1}y^{1-l}$  for all $i=2,3,4$ (and $n=5$).
\end{description}

Consider the possibilities (i)--(iv) for $g_2$. Since no two of $g_2$, $g_2^{-1}x$ and $g_2^{-1}y$ can
be conjugate in $K$ we must have (i) and (iii). 
Then $g_2= (x^{-k} \delta^{-m} y^l)\, g_2\, (y^{-l} \delta^m x^k)$. 
The element $y^{-l} \delta^m x^k$ of $K$ commutes with $g_2$, thus it must be a power of $g_2$. That
is to say $\delta^m = y^l g_2^q x^{-k}$ for some $q\in\Z$. 
In particular, since $\delta^2,x,y,g_1\in F_{n-1}$, but $\delta\notin F_{n-1}$, we must have $m$ even. 
Now, take $x, g_2, g_3, \dots, g_{n-2}, y$ as a free basis for $F_{n-1}$. With respect to these generators, 
$z=\delta^2$ has a reduced form
\[
z= x\, y^{-1}\,g_2 g_4 \dots g_{n-3}\,(g_2g_3  \dots g_{n-1} g_{n-2})^{-1}\, g_3 g_5 \dots g_{n-2}\ x^{-1}\, y\,.
\]
Therefore, $\delta^m=z^{m \over 2}$ cannot have the form $y^l g_2^q x^{-k}$ unless $k=l=q=m=0$. So, 
$\varphi$ must be the identity on $F_{n-1}$. This completes the proof of Theorem \ref{autos}.
\endproof

We shall now work with the presentation $A(D_n)=F_{n-1}\rtimes{\rho_D}\Bn$ of Section \ref{sect1}.
Note that $A(D_n)$ is an index 2 normal subgroup of $K\rtimes_{\rho^+}\Bn\cong\pi_1N^+$,
where $N^+$ was introduced in Section \ref{sect2}. (In fact $K\rtimes\Bn\cong A(D_n)\rtimes C_2$ where
the section sends the generator of $C_2$ to $x_1$). 
For notational convenience, we shall identify elements of $\Bn<A(D_n)$
with the corresponding elements of $\Inn(A(D_n))$ via their action by conjugation 
(and using the fact that $\rho_D$ is faithful). We call these elements \emph{braid automorphisms}.
As with the $B_n$ case, it is easily seen that $\zeta$ acts on $K\rtimes\Bn$, and so on $A(D_n)$,
 by conjugation by $\delta^{-1}$, where $\de =x_1x_2...x_n$.
It is also a straightforward exercise to check that the centre of $K\rtimes\Bn$ is generated by the
element $\de\zeta$, and that the centre of $A(D_n)$ is generated by $\de\zeta$ if $n$ is even, and 
$\de^2\zeta^2$ if $n$ is odd. 

\paragraph{Special automorphisms.}
We define the automorphism $\ep_n\in\Aut(K\rtimes\Bn)$ such that $\ep_n(\al_i)=\al_i^{-1}$
for $i=1,2,..,n-1$, and $\ep_n(x_1)=x_1$. It is a short exercise to show that $\ep_n$ is well-defined.
Moreover, $\ep_n$ induces the automorphism of $A(D_n)$ 
which, by abuse of notation, we shall also call $\ep_n$, and which is defined by
$\ep_n(\de_i)=\de_i^{-1}$ for every $i=1,2,..,n$.

We also define the so-called \emph{graph automorphism} $\tau_n\in\Aut(A(D_n))$
induced by the involution of the Coxeter graph of type $D_n$, namely 
$\tau_n(\de_1)=\de_2$, $\tau_n(\de_2)=\de_1$, and $\tau_n(\de_i)=\de_i$ if $i\geq 3$.

Note that, since $A(D_n)$ is a normal subgroup of $K\rtimes\Bn$,
every inner automorphism of $K\rtimes\Bn$  induces an automorphism of $A(D_n)$.
These constitute a subgroup of $\Aut(A(D_n))$ which is generated by the inner
automorphisms of $A(D_n)$  and conjugation by $x_1$ in the larger group.
The latter automorphism of $A(D_n)$ is precisely the graph automorphism
$\tau_n$.

Finally, we note that the automorphism $\tau_n$ is an inner automorphism of $A(D_n)$
if and only if $n$ is odd.  This is a straightforward consequence of \cite{Par}. Alternatively,
one can observe that if $\tau_n$ is an inner automorphism, i.e: conjugation by some element $k\be$, say,
where $k\in F_{n-1}$ and $\be\in\Bn$, then the element $x_1k\be\in K\rtimes\Bn$ centralizes $A(D_n)$.
But that is to say that conjugation by $x_1k$ in $K$ agrees with the action of $\be^{-1}$, and so 
fixes $\de$. But then $x_1k\in K\setminus F_{n-1}$ is a power of $\de$, which is only possible
if $n$ is odd. On the other hand, if $n$ is odd, $x_1$ differs from $\de$ by an element
of $F_{n-1}$, and conjugation by $\de$ is realised by the braid automorphism $\zeta$,
so $\tau_n$ is an inner automorphism of $A(D_n)$.

\begin{thm}\label{AutKFDn}
Let $n\geq 4$.
The group $\Aut(A(D_n),F_{n-1})$ of automorphisms of $A(D_n)$ leaving invariant
the subgroup $F_{n-1}$ is generated by the inner automorphisms, the graph automorphism $\tau_n$,
and the automorphism $\ep_n$. More precisely ,
\[
\Aut (A(D_n), F_{n-1}) = 
\begin{cases}
\Inn (A(D_n)) \rtimes \< \epsilon_n, \tau_n \> \cong \big( A(D_n)/ Z(A(D_n)) \big) \rtimes (C_2 
\times C_2) &\text{ if } n \text{ is even}\\
\Inn (A(D_n)) \rtimes \< \epsilon_n \> \cong \big( A(D_n)/ Z(A(D_n)) \big) \rtimes C_2 
&\text{ if } n \text{ is odd}
\end{cases}
\]
where $Z(A(D_n))$ denotes the centre of $A(D_n)$. In particular, 
$\Out(A(D_n),F_{n-1}) \cong C_2 \times C_2$ if $n$ is even,
and $\Out(A(D_n),F_{n-1}) \cong C_2$ if $n$ is odd.
\end{thm}

\Proof
Let $\varphi$ be an arbitrary automorphism of $A(D_n)=F_{n-1}\rtimes\Bn$
such that $\varphi(F_{n-1})=F_{n-1}$. 
Then $\varphi$ induces an automorphism $\overline{\varphi}$ of the quotient group
$\Bn=A(D_n)/F_{n-1}$. Using the theorem of Dyer and Grossman \cite{DG}, we may suppose,
up to a multiplication of $\varphi$ by a braid automorphism, and by $\ep_n$ if necessary,
that $\overline{\varphi}$ is the trivial automorphism of $\Bn$. In other words, there exists
some function $k:\Bn\to F_{n-1}$, written $\ga\mapsto k_\ga$,
such that
\[
\varphi(\ga)= k_\ga\,\ga\,,\hskip5mm \text{ for all } \ga\in\Bn\,.
\] 
We now consider the element $z\in K$ which generates the subgroup of fixed
points in $F_{n-1}$ under the braid action. That is $z=\de$, if $n$ is even, and $\de^2$
if $n$ is odd. For any $\ga\in\Bn$ we have
\begin{equation}\label{phideleqn} 
\varphi(z)=\varphi(\ga z \ga^{-1})=k_\ga\,\ga\,\varphi(z)\,\ga^{-1}k_\ga^{-1}\,.
\end{equation}
That is to say, $\ga(\varphi(z))\sim\varphi(z)$ for all $\ga\in\Bn$, where $\sim$ denotes conjugacy in $K$.
But then Lemma \ref{deltakey}, plus the fact that $z$ is not a proper power in $F_{n-1}$,
implies that $\varphi(z)$ is conjugate in $K$ to $z^\nu$, for $\nu\in\{\pm 1\}$.

Now, multiplying $\varphi$ by an inner automorphism from $F_{n-1}$ and by $\tau_n$ if necessary 
(that is by some automorphism which is induced by conjugation, in the bigger group $K\rtimes\Bn$,
by an element of $K$), we may suppose, in fact, that $\varphi(z)=z^\nu$.
Note that this modification of $\varphi$ does not change the induced
automorphism $\overline{\varphi}$, but it does modify the function $k$.
In any case, the equation (\ref{phideleqn}) together with the fact that braid automorphisms fix $z$,
gives simply $z^\nu=k_\ga z^\nu k_\ga^{-1}$. But then $k_\ga$ must be a power of $z$ for 
every $\ga\in\Bn$. In fact, $k:\Bn\to\< z\>\cong \Z$ must be a homomorphism, since 
$\varphi(\ga\be)=k_\ga\,\ga\, k_\be\,\be=k_\ga k_\be\,\ga\be$.
Moreover, in view of the braid relations $\al_i\al_{i+1}\al_i=\al_{i+1}\al_i\al_{i+1}$,
the only homomorphisms $\Bn\to\Z$ are multiples of the length homomorphism $\ell:\Bn\to\Z$
which sends each standard generator to 1. In other words, there exists an $m\in\Z$ such that
$\varphi(\ga)=z^{m\ell(\ga)}\ga$, for all $\ga\in\Bn$. 

Finally, we observe that 
the centre of $A(D_n)$ is generated by the element $z\zeta$ if $n$ is even, and 
$z\zeta^2$ if $n$ is odd, and is left invariant by
any automorphism of $A(D_n)$. For instance, if $n\geq 5$ is odd,
$\varphi(z\zeta)=(z\zeta)^{\pm 1}$. 
Since $\varphi(z\zeta)=z^\nu\, z^{m\ell(\zeta)}\,\zeta$ we deduce that 
$m\ell(\zeta)=1-\nu$. But, since $\ell(\zeta)=n(n-1)\geq 20$, we must have $m=0$. Similarly,
when $n$ is even we deduce that $m=0$, and in both cases the action of $\varphi$ on $F_{n-1}$
is therefore $\Bn$-equivariant. By Theorem \ref{autos}, we now have that $\varphi$ agrees on $F_{n-1}$
with a central braid automorphism $\zeta^k$ for some $k\in\Z$. But then, $\zeta^{-k}\varphi$
fixes both subgroups $F_{n-1}$ and $\Bn$, so is trivial, and $\varphi$ is an inner automorphism.
\endproof

\subsection{Automorphisms of $K \rtimes \Bn$}

We continue with the same notations as in the previous subsection.
Recall that $K=\star_n (C_2)$ denotes the free product of $n$ copies of $C_2$, and $\Bn$ 
acts on $K$ via the representation $\rho^+$ define in Section \ref{sect1}. We proceed now to 
determine the automorphism group of $K \rtimes_{\rho^+} \Bn$ following the same 
strategy as in the previous subsections. First, we establish the following result which is 
analogous to Proposition \ref{Bequivautos} for the group $A(B_n)$ and to Theorem \ref{autos} for 
the group $A(D_n)$.

\begin{prop}\label{an1}
Let $n \ge 3$, and identify $\Bn$ with a subgroup of $\Aut(K)$ via the representation $\rho^+$. 
Then the group of $\Bn$-equivariant automorphisms of $K$ is the cyclic subgroup of $\Aut(K)$ 
generated by $\zeta$.
\end{prop}

\Proof
Let $\varphi: K \to K$ be a $\Bn$-equivariant automorphism. By Lemma \ref{deltakey}, $K^\Bn$ is 
the infinite cyclic subgroup generated by $\delta = x_1 x_2 \dots x_n$, thus $\varphi( \delta)= 
\delta^\nu$ where $\nu \in \{ \pm 1 \}$. Let $T= \{\alpha_2, \alpha_3, \dots, \alpha_{n-1} \}$. 
By Lemma \ref{second}, we have $K^{\<T\>} = \< x_1, x_2 x_3 \dots x_n \> = \<x_1\> \star 
\<\delta\> \cong C_2 \star \Z$, thus $\varphi$ restricts to an automorphism $\<x_1\> \star 
\<\delta\> \to \<x_1\> \star \<\delta\>$. The element $\varphi(x_1)$ is of order $2$ and 
$\varphi(x_1) \in K^{\<T\>}$, thus $\varphi(x_1)$ is conjugate to $x_1$ in $K^{\<T\>}$, namely, 
there exists $w \in K^{\<T\>}$ such that $\varphi(x_1) = w x_1 w^{-1}$. Now, $\varphi(\delta) = 
\delta^{\pm 1}$ and $\varphi(x_1)= wx_1w^{-1}$ generate $K^{\<T\>}= \<x_1\> \star \<\delta\>$, 
and this is possible only if $w$ is of the form $w= \delta^k x_1^\mu$ with $k \in \Z$ and $\mu\in 
\{0,1\}$, thus $\varphi(x_1)= \delta^k x_1 \delta^{-k}$ for some $k \in \Z$. Observe that, for 
$i=1,2, \dots, n$, there exists a braid $\gamma_i \in \Bn$ such that $\gamma_i (x_1) = x_i$. By 
equivariance, it follows that $\varphi(x_i)= \delta^k x_i \delta^{-k}$ for all $i=1, \dots, n$, 
and therefore $\varphi= \zeta^{-k}$.
\endproof

Recall that $\epsilon_n$ denotes the automorphism of $K \rtimes \Bn$ determined by 
$\epsilon_n(x_1)=x_1$ and $\epsilon_n( \alpha_i) = \alpha_i^{-1}$ for $i=1, \dots, n-1$.

\begin{thm}\label{an2}
Let $n \ge 2$. The group $\Aut( K \rtimes \Bn)$ is generated by the inner automorphisms and the 
automorphism $\epsilon_n$. Thus 
$$
\Aut (K \rtimes \Bn) = \Inn (K \rtimes \Bn)\rtimes \< \epsilon_n\> \cong \big( (K \rtimes \Bn)/ Z(K \rtimes \Bn) \big) 
\rtimes C_2\,,
$$
where $Z(K \rtimes \Bn)$ denotes the centre of $K \rtimes \Bn$. In particular, $\Out(K \rtimes \Bn) \cong C_2$.
\end{thm}

\Proof
The case $n=2$ is special and should be treated separately. We leave this case to the reader, and assume
from now on that $n\ge 3$.

Let $\varphi: K \rtimes \Bn \to K \rtimes \Bn$ be an automorphism. Note that $K$ is the smallest 
subgroup of $K \rtimes \Bn$ which contains every element of order $2$, and therefore is a 
characteristic subgroup. Let $\ov{\varphi}$ be the automorphism of $\Bn$ induced by $\varphi$. 
Using the theorem of Dyer and Grossman \cite{DG}, we may suppose, up to a multiplication of 
$\varphi$ by a braid automorphism and by $\epsilon_n$ if necessary, that $\ov{\varphi}$ is the 
trivial automorphism of $\Bn$. So, there exists a function $k: \Bn \to K$, $\gamma \mapsto 
k_\gamma$, such that $\varphi (\gamma) = k_\gamma \gamma$ for all $\gamma \in \Bn$.

Recall that every braid fixes $\delta=x_1x_2 \dots x_n$. For any $\gamma \in \Bn$, we have
\begin{equation}\label{an3}
\varphi (\delta) = \varphi (\gamma \delta \gamma^{-1}) = k_\gamma \gamma \varphi(\delta) 
\gamma^{-1} k_\gamma^{-1}\,.
\end{equation}
This shows that $\gamma (\varphi (\delta)) \sim \varphi(\delta)$ for all $\gamma \in \Bn$, thus, 
by Lemma \ref{deltakey}, $\varphi(\delta)$ is conjugate to $\delta^\nu$ for some $\nu \in \{\pm 
1\}$. Multiplying $\varphi$ by an inner automorphism from $K$, we may suppose that 
$\varphi(\delta) = \delta^\nu$. This modification does not change the induced automorphism 
$\ov{\varphi}$, but it does modify the function $k: \Bn \to K$.

The equation (\ref{an3}) together with the fact that braid automorphisms fix $\delta$, gives 
$\delta^\nu= k_\gamma \delta^\nu k_\gamma^{-1}$, thus $k_\gamma$ is a power of $\delta$. It 
follows that the function $k: \Bn \to \<\delta\> \cong \Z$ is a homomorphism since, for $\gamma, 
\beta \in \Bn$, $k_\beta$ being a power of $\delta$, we have $\varphi (\gamma \beta) = k_\gamma 
\gamma k_\beta \beta= k_\gamma k_\beta \gamma \beta$. The only homomorphisms $\Bn \to \Z$ are 
multiples of the length function $\ell: \Bn \to \Z$, thus there exists $m \in \Z$ such that 
$k_\gamma = \delta^{m \ell (\gamma)}$ for all $\gamma \in \Bn$.

Recall that the centre of $K \rtimes \Bn$ is the infinite cyclic subgroup generated by $\delta 
\zeta$. In particular, we have $\varphi (\delta \zeta) = (\delta \zeta)^{\pm 1}$. Since $\varphi 
(\delta \zeta) = \delta^\nu \delta^{m \ell (\zeta)} \zeta$, it follows that $m \ell (\zeta)= 1-
\nu$, thus $m=0$ (since $\ell (\zeta) \ge 6$). This shows that $\varphi(\gamma)=\gamma$ for all 
$\gamma \in \Bn$, and therefore that the action of $\varphi$ on $K$ is $\Bn$-equivariant. By 
Proposition \ref{an1}, $\varphi$ agrees on $K$ with a central braid automorphism $\zeta^k$ for 
some $k \in \Z$. But then, $\zeta^{-k} \varphi$ fixes both subgroups $K$ and $\Bn$, so it is 
trivial, and $\varphi$ is an inner automorphism.
\endproof


\section{Automorphisms of Artin groups of rank 2} \label{chap5}

Let $m$ be a positive integer, $m \ge 3$. The aim of this section is to determine the group of 
automorphisms of the Artin group
$$
A= \< \alpha, \beta \ |\ w(\alpha, \beta:m) = w(\beta, \alpha :m) \>\,.
$$
We start by describing some automorphisms.

\paragraph{Special automorphisms.}
We define the automorphism $\epsilon \in \Aut(A)$ by $\epsilon (\alpha) = \alpha^{-1}$ and 
$\epsilon (\beta) = \beta^{-1}$. We also define the so-called {\it graph automorphism} $\tau \in 
\Aut(A)$ by $\tau(\alpha) = \beta$ and $\tau (\beta) = \alpha$. Note that $\tau$ is inner if and 
only if $m$ is odd, and $\tau$ commutes with $\epsilon$. If $m$ is odd, then $\tau$ is simply the 
conjugation by the element $\Delta= w(\alpha, \beta :m) \in A$.

Assume $m$ is even. It is an easy exercise to show that the mapping $\alpha \mapsto \beta^{-1}$, 
$\beta \mapsto \beta \alpha \beta$, determines an automorphism $\eta: A \to A$ of infinite order. 
Moreover, $\eta$ and $\epsilon$ commute, and $\tau \eta \tau = \eta^{-1}$.

\begin{thm}\label{thm51}
\begin{description}
\item{(i)} Let $m$ be odd. Then $\Aut(A)$ is generated by the inner automorphisms and the automorphism 
$\epsilon$. Thus
$$
\Aut(A) = \Inn(A) \rtimes \< \epsilon \> \cong (A/Z(A)) \rtimes C_2\,,
$$
and, in particular, $\Out(A) \cong C_2$.
\item{(ii)} Let $m$ be even. Then $\Aut(A)$ is generated by the inner automorphisms and the automorphisms 
$\epsilon, \tau, \eta$. Thus
$$
\Aut(A) = \Inn(A) \rtimes \< \epsilon, \tau, \eta \> \cong (A/Z(A)) \rtimes ((\Z \rtimes C_2) 
\times C_2)\,,
$$
and,  in particular, $\Out(A) \cong (\Z \rtimes C_2) \times C_2$.
\end{description}
\end{thm}

In order to prove Theorem \ref{thm51}, we first need to determine the automorphisms of the groups 
$C_2 \star C_m$ and $C_k \star \Z$ ($k \ge 2$). These are laid out in the following two Lemmas.

Let $G$ be a group. For $g \in G$, we denote by $\iota_g: G \to G$ the inner automorphism which 
sends $h$ to $ghg^{-1}$ for all $h \in G$. Let $\varphi_1, \varphi_2 : G \to G$ be two 
automorphisms. We say that $\varphi_1$ and $\varphi_2$ are {\it conjugate} and denote $\varphi_1 
\sim \varphi_2$ if there exists $g \in G$ such that $\varphi_2 = \varphi_1 \iota_g$.

Let $u$ and $v$ be the standard generators of $C_2 \star C_m$. That is, $C_2 \star C_m = \<u,v| 
u^2=v^m=1\>$. For any $r \in \{ 1,2, \dots, m-1\}$ such that $r$ and $m$ are coprime, we define the 
automorphism $\nu_r: C_2 \star C_m \to C_2 \star C_m$ by $\nu_r (u)=u$ and $\nu_r(v) =v^r$.

\begin{lemma}\label{lemma52}
Let $m\ge 3$.
Let $\varphi$ be an automorphism of $C_2 \star C_m$. Then there exists $r \in \{1,2, \dots, m-
1\}$ such that $r$ and $m$ are coprime and $\varphi \sim \nu_r$.
\end{lemma}

\Proof
The element $\varphi(u)$ is of order two, thus it is conjugate to $u$. So, up to multiplication 
of $\varphi$ by an inner automorphism, we may suppose that $\varphi(u) =u$. On the other hand, 
$\varphi(v)$ is of order $m \ge 3$, thus it is conjugate to some $v^r$, where $r\in \{1,2, \dots, 
m-1\}$ and $r$ and $m$ are coprime. Write $\varphi(v)= w v^r w^{-1}$, where $w\in C_2 \star C_m$. 
Now the inclusion $v \in \im \varphi$ is possible only if $w$ is of the form $w=u^\varepsilon 
v^p$, where $\varepsilon \in \{0,1\}$ and $p \in \Z$. So, $\varphi(v)= u^\varepsilon v^r u^{-
\varepsilon}$, and therefore $\varphi= \iota_u^\varepsilon\, \nu_r=\nu_r\, \iota_u^\varepsilon$.
\endproof

Let $u,v$ be the standard generators of $C_k \star \Z$. That is, $C_k \star \Z = \< u,v| u^k=1 
\>$. Choose $\varepsilon\in \{\pm 1\}$ and $r,s \in \{0,1, \dots, k-1\}$ such that $r$ and $k$ 
are coprime, and define the automorphism $\nu_{\varepsilon, r, s}: C_k \star \Z \to C_k \star 
\Z$ by $\nu_{\varepsilon, r, s} (u) = u^r$ and $\nu_{\varepsilon, r, s} (v) = v^{\varepsilon} 
u^s$.

\begin{lemma}\label{lemma53}
Let $k \ge 2$.
Let $\varphi$ be an automorphism of $C_k \star \Z$. Then there exist $\varepsilon \in \{ \pm 1\}$ 
and $r,s \in \{ 0,1, \dots, k-1\}$ such that $r$ and $k$ are coprime and $\varphi \sim 
\nu_{\varepsilon, r, s}$.
\end{lemma}

\Proof
The element $\varphi(u)$ is of order $k$, thus it is conjugate to some $u^r$ where $r \in \{1,2, 
\dots, k-1\}$ and $r$ and $k$ are coprime. So, up to multiplication of $\varphi$ by an inner 
automorphism, we may assume that $\varphi(u)=u^r$.
The element $\varphi(v)$ is of infinite order. It can be written
$$
\varphi (v) = u^{n_0} v^{m_1} u^{n_1} \dots v^{m_\ell} u^{n_\ell}\,,
$$
where $0\le n_0,n_\ell \le k-1$, $1\le n_i\le k-1$ for $i=1, \dots, \ell-1$, $m_i \in \Z\setminus \{0\}$ for 
$i=1, \dots, \ell$, and $\ell \ge 1$. Now, by exactly the same argument as in the proof of Lemma 
\ref{third}, one easily verifies that $v$ is in the subgroup 
generated by $\varphi(u)=u^r$ and $\varphi(v)$ only if $\ell=1$ and $m_1\in \{ \pm 1\}$. It follows 
that $\varphi= \iota_u^{n_0}\, \nu_{\varepsilon, r, s}\,$, where $\varepsilon=m_1$ and 
$s\equiv n_0+n_1\ (k)$.
\endproof

\paragraph{Proof of Theorem \ref{thm51}, case $m$ odd.}
Write $m=2k+1$. Then
$$
A=\< \alpha, \beta\ |\ (\alpha \beta)^k \alpha = \beta (\alpha \beta)^k\>\,.
$$
Let $a=\alpha \beta$ and $b= (\alpha \beta)^k \alpha = \beta (\alpha \beta)^k$. The elements $a$ 
and $b$ generate $A$ and a presentation for $A$ with respect to these generators is
$$
A=\< a,b\ |\ a^m=b^2\>\,.
$$
Note also that $\tau \epsilon$ is the automorphism of $A$ which sends $a$ to $a^{-1}$ and $b$ to 
$b^{-1}$.
Let $c=a^m=b^2$. A straightforward consequence of the above presentation is that the centre of 
$A$, denoted by $Z(A)$, is the infinite cyclic subgroup of $A$ generated by $c$. Let $\ov{A}= 
A/Z(A)$. Then $\ov{A}=\<\ov{a}, \ov{b}| \ov{a}^m = \ov{b}^2 =1 \> \cong C_m \star C_2$.

Let $\varphi: A \to A$ be an automorphism. Let $\ov{\varphi}$ be the automorphism of $\ov{A}$ 
induced by $\varphi$. By Lemma \ref{lemma52}, up to multiplication of $\varphi$ by an inner 
automorphism if necessary, we may assume that $\ov{\varphi} (\ov{a}) = \ov{a}^r$ and $\ov{\varphi} 
(\ov{b}) = \ov{b}$ where $r \in \{1,2, \dots, m-1\}$ (and where $r$ and $m$ are coprime, although we
will not need to use this fact). In other words, there exists $p,q \in \Z$ such that
$$
\varphi(a) = a^rc^p = a^{r+mp}\ \ \text{ and }\ \ \varphi(b) = bc^q =b^{1+2q}\,.
$$
Recall that the centre of $A$ is the infinite cyclic subgroup generated by $c$. In particular, 
$\varphi(c)=c$ or $c^{-1}$. Replacing $\varphi$ by $(\tau\epsilon)\varphi$ if necessary we may assume
that $\varphi(c)=c$. Then the equality $c=\varphi(c)=\varphi(a)^m= a^{m(r+mp)} = 
c^{r+mp}$ implies that $r+mp=1$, and therefore $\varphi(a)=a$.
Similarly we observe that $\varphi(b)=b$, and so $\varphi$ must be the trivial automorphism of $A$.
\endproof

\paragraph{Proof of Theorem \ref{thm51}, case $m$ even.}
Write $m=2k$. Then
$$
A=\< \alpha, \beta \ |\ (\alpha \beta)^k= (\beta \alpha)^k \>\,.
$$
Let $a=\alpha \beta$ and $b=\beta$. Then the elements $a$ and $b$ generate $A$ and a presentation 
for $A$ with respect to these generators is
$$
A=\< a,b\ |\ a^k b=b a^k \>\,.
$$
Note also that $(\eta \tau \epsilon)$ is the automorphism of $A$ which sends $a$ to $a^{-1}$ and 
$b$ to $b$, $(\iota_\beta^{-1} \eta \tau)$ is the automorphism which sends $a$ to $a$ and $b$ to $b^{-
1}$, and $\eta$ is the automorphism which sends $a$ to $a$ and $b$ to $ba$.
Let $c=a^k$. It follows from the above presentation that the centre of $A$ is the infinite cyclic 
subgroup generated by $c$. Let $\ov{A} = A/Z(A)$. Then $\ov{A}= \< \ov{a}, \ov{b}| \ov{a}^k=1\> 
\cong C_k \star \Z$.

Let $\varphi: A \to A$ be an automorphism. Let $\ov{\varphi}$ be the automorphism of $\ov{A}$ 
induced by $\varphi$. By Lemma \ref{lemma53}, up to multiplication of $\varphi$ by an inner 
automorphism, we may assume that $\ov{\varphi} (\ov{a})= \ov{a}^r$ and $\ov{\varphi} (\ov{b}) = 
\ov{b}^\varepsilon \ov{a}^s$, where $\varepsilon \in \{ \pm 1\}$, $r,s \in \{ 0,1, \dots, k-1\}$, 
and $r$ and $k$ are coprime. It follows that there exist $p,q \in \Z$ such that
$$
\varphi(a)= a^r c^p = a^{r+kp}\,, \qquad \varphi(b)= b^\varepsilon a^s c^q = b^\varepsilon 
a^{s+kq}\,.
$$
Recall that the centre of $A$ is generated by $c=a^k$. Thus $\varphi(c)=c$ or $c^{-1}$. 
Replacing $\varphi$ by $(\eta \tau \epsilon)\varphi$ if necessary, we may assume that $\varphi(c)=c$. 
In this case, the equality $c=\varphi(c) = \varphi(a)^k =c^{r+kp}$ implies that $r+kp=1$, 
so that $\varphi(a)=a$. Now replacing $\varphi$ by  $(\iota^{-1}_\beta \eta \tau) \varphi$ if necessary,
we may assume $\varep=1$. But then $\varphi(a)=a$ and $\varphi(b)=ba^{s+kq}$, so $\varphi= \eta^{s+kq}$.
\endproof


\end{document}